\documentclass{article}

\usepackage{pbmath}
\usepackage{epsf}
\usepackage{psfrag}
\usepackage{amsmath}
\usepackage{amsthm}
\usepackage{amssymb}

\renewcommand{\epsfsize}[2]{0.7\textwidth}

\title{Dynamics of free group automorphisms}
\date{\today}
\author{Peter Brinkmann}

\begin{document}
\maketitle

% $Id: abstract.tex,v 1.3 2008/06/01 05:04:29 brinkman Exp brinkman $
\begin{abstract}
We present a coarse convexity result for the dynamics of free group
automorphisms: Given an automorphism $\phi$ of a finitely generated free
group $F$, we show that for all $x\in F$ and $0\leq i\leq N$, the length
of $\phi^i(x)$ is bounded above by a constant multiple of the sum of the
lengths of $x$ and $\phi^N(x)$, with the constant depending only on
$\phi$.
\end{abstract}

\section*{Introduction}

The following theorem is the main result of this paper. It follows from
a technical result (\thmref{mainthm}) that uses the machinery of
improved relative train track maps of Bestvina, Feighn, and Handel
\cite{tits1}.

\begin{thm}\label{mainapp}
Let $\phi\co F \rightarrow F$ be an automorphism of a finitely generated
free group. Then there exists a constant $K\geq 1$ such that for any
pair of exponents $N,i$ satisfying $0\leq i \leq N$, the following two
statements hold:

\begin{enumerate}
\item If $w$ is a cyclic word in $G$, then
\[
	||\phi^i_\#(w)||\leq K\left(||w||+||\phi^N_\#(w)||\right),
\]
where $||w||$ is the length of the cyclic reduction of $w$ with
respect to some word metric on $F$.
\item If $w$ is a word in $F$, then
\[
	|\phi^i_\#(w)|\leq K\left(|w|+|\phi^N(w)|\right),
\]
where $|w|$ is the length of $w$.
\end{enumerate}
Given an improved relative train track representative of some power of
$\phi$, the constant $K$ can be computed.
\end{thm}

\begin{rem}[A note on computability]\label{compnote}
Given an automorphism $\phi\co F\rightarrow F,$ we can compute a
relative train track representative of $\phi$ \cite{hb1, dv}. The
construction of {\em improved} relative train track maps, however,
involves a compactness argument in a universal cover \cite[Proof of
Proposition~5.4.3]{tits1} that is not constructive.  A number of
algorithmic improvements of relative train tracks appear in
\cite{pborbits}, in the context of an algorithm that detects automorphic
orbits in free groups.
\end{rem}

The statement of the theorem does not depend on the choice of generators
of $F$.  The intuitive meaning of the theorem is that the map $i\mapsto
|\phi^i(w)|$ is coarsely convex for all words $w\in F$.  Klaus Johannson
informed me that a similar result is a folk theorem in the case of
surface homeomorphisms. Also, while free-by-cyclic groups are not, in
general, ${\rm CAT}(0)$-groups \cite{notcat0}, \thmref{mainapp} suggests
that their dynamics mimics that of ${\rm CAT}(0)$-groups.
\thmref{mainapp} complements the following strong convexity result in an
important special case.

\begin{thm}[\cite{pbhyp}]
If $\phi\co F\rightarrow F$ is an {\em atoroidal} automorphism, i.e.,
$\phi$ has no nontrivial periodic conjugacy classes, then $\phi$ is {\em
hyperbolic}, i.e., there exists a constant $\lambda>1$ such that
\[ |x|\leq \lambda\max\left\{|\phi^{\pm 1}(x)|\right\} \]
for all $x\in F$.
\end{thm}

I originally set out to prove \thmref{mainapp}
because it immediately implies that in a free-by-cyclic group
\[
    \Gamma=F \rtimes_\phi \mathbb{Z}=
    \langle~x_1,\ldots,x_n,t~|~t^{-1}x_it=\phi(x_i)~\rangle,
\]
words of the form $t^{-k}wt^k\phi^k(w^{-1})$ satisfy a quadratic
isoperimetric inequality. (Note, however, that \thmref{mainapp} is
stronger than the mere existence of a quadratic isoperimetric inequality
for such words.) Natasa Macura previously proved a quadratic
isoperimetric inequality for mapping tori of automorphisms of
polynomial growth \cite{natasagafa}.  Martin Bridson and Daniel Groves
have since proved that all free-by-cyclic groups satisfy a quadratic
isoperimetric inequality \cite{bgip}.  They also obtain a new proof
of \thmref{mainapp} as an application of their techniques.

In \secref{train}, we review the pertinent definitions and results
regarding train track maps from \cite{tits1}. We also state the main
technical result, \thmref{mainthm}, and we show how \thmref{mainapp}
follows from \thmref{mainthm}.  \secref{train2} provides some more
results on train tracks and automorphisms of free groups. \secref{exsec}
introduces some notation and terminology and lists a number of examples
that illustrate some of the issues and subtleties that need to be
addressed in the proof of \thmref{mainthm}.  \secref{nonlinsec}
establishes a technical proposition that may be of independent interest.
Finally, \secref{polysec} and \secref{pfsec} contain the proof of
\thmref{mainthm}.

I would like to express my gratitude to Ilya Kapovich for many helpful
discussions, to Mladen Bestvina for patiently answering my questions,
to Steve Gersten for encouraging me to write up this result for its
own sake, to the University of Osnabr\"uck for their hospitality,
and to Swarup Gadde and the University of Melbourne as well as the
Max-Planck-Institute of Mathematics for their hospitality and
financial support. Klaus Johannson and Richard Weidmann kindly served
as a sounding board while I was working on the exposition of this paper.

\section{Improved relative train track maps}\label{train}

In this section, we review the theory of train tracks developed in
\cite{hb1,tits1}. We will restrict our attention to the collection of
those results that we will use in this paper.

Given an automorphism $\phi\in Aut(F)$, we can find a based homotopy
equivalence $f\co G\rightarrow G$ of a finite connected graph $G$ such
that $\pi_1(G)\cong F$ and $f$ induces $\phi$. This observation allows us to
apply topological techniques to automorphisms of free groups. In many
cases, it is convenient to work with outer automorphisms. Topologically,
this means that we work with homotopy equivalences rather that based
homotopy equivalences.

Oftentimes, a homotopy equivalence $f\co G\rightarrow G$ will respect a
{\em filtration}
of $G$, i.\ e., there exist subgraphs $G_0=\emptyset\subset G_1 \subset \cdots
\subset G_k=G$ such that for each filtration element $G_r$,
the restriction of $f$ to $G_r$
is a homotopy equivalence of $G_r$. The subgraph
$H_r=\overline{G_r\setminus G_{r-1}}$ 
is called the {\em $r$-th stratum} of the filtration.
We say that a path $\rho$ has {\em nontrivial intersection} with a
stratum $H_r$ if $\rho$ crosses at least one edge in $H_r$.

If $E_1,\cdots,E_m$ is the collection of edges in some stratum $H_r$, the
{\em transition matrix} of $H_r$ is the nonnegative
$m\times m$-matrix $M_r$ whose $ij$-th
entry is the number of times the $f$-image of $E_j$ crosses $E_i$, regardless
of orientation. $M_r$ is said to be {\em irreducible} if
for every tuple $1\leq i,j \leq m$, there exists some exponent $n>0$ such that
the $ij$-th entry of $M_r^n$ is nonzero.
If $M_r$ is irreducible, then it has a maximal real eigenvalue
$\lambda_r\geq 1$
\cite{gantmacher}.  We call $\lambda_r$ the
{\em growth rate} of $H_r$.

Given a homotopy equivalence $f\co G\rightarrow G$, we can always find a
filtration of $G$ such that each transition matrix is either a zero matrix
or irreducible. A stratum $H_r$ in such a filtration is called
{\em zero stratum} if $M_r=0$. $H_r$ is called {\em exponentially growing}
if $M_r$ is irreducible with $\lambda_r>1$, and it is called
{\em polynomially growing} if $M_r$ is irreducible with $\lambda_r=1$.

An unordered pair of edges in $G$ originating from the same vertex is called
a {\em turn}. A turn is called {\em degenerate} if the two edges are equal.
We define a map
$Df\co \{\text{turns in } G\}\rightarrow \{\text{turns in } G\}$ by sending each
edge in a turn to the first edge in its image under $f$. A turn is called
{\em illegal} if its image under some iterate of $Df$ is degenerate, {\em legal}
otherwise.

An edge path $\rho=E_1E_2\cdots E_s$ is said to contain the turns
$(E_i^{-1},E_{i+1})$ for $1\leq i <s$. $\rho$ is said to be legal if
all its turns are legal, and a path $\rho\subset G_r$ is $r$-legal if
no illegal turn in $\alpha$ involves an edge in $H_r$.

Let $\rho$ be a path in $G$. In general, the composition
$f^k\circ \rho$ is not an immersion, but there is exactly one immersion
that is homotopic to $f^k\circ \rho$ relative endpoints. We denote this
immersion by $f^k_\#(\rho)$, and we say that we obtain $f^k_\#(\rho)$ from
$f^k\circ \rho$ by {\em tightening}. If $\sigma$ is a circuit in $G$,
then $f^k_\#(\sigma)$ is the immersed circuit homotopic to $f^k\circ \sigma$.

\begin{rem}\label{ambiguity}
A path is tightened by cancelling adjacent pairs of inverse edges until
no inverse pairs are left. The result of such a sequence of
cancellations is uniquely determined, but the sequence is not. For
instance, $EE^{-1}E$ may be tightened as $E(E^{-1}E)$ or $(EE^{-1})E$.
\end{rem}

\begin{conv}\label{howtotighten}
Let $\rho_i, i=1,\ldots,k$ be paths that can be concatenated to form a
path $\rho=\rho_1\rho_2\cdots\rho_k$. When tightening $f(\rho)$ to obtain
$f_\#(\rho)$, we adopt the convention that we first tighten the images of
$\rho_i$ to $f_\#(\rho_i).$ In a second step, we tighten the concatenation
$f_\#(\rho_1)\cdots f_\#(\rho_k)$ to $f_\#(\rho).$

In many situations, the length of a subpath $\rho_i$ will be greater
than the number of edges that cancel at either end, in which case it
makes sense to talk about edges in $f_\#(\rho)$ {\em originating from
$\rho_i.$}
\end{conv}

A path $\rho$ is a {\em (periodic) Nielsen path} if $f^k_\#(\rho)=\rho$
for some $k>0$. In this case, the smallest such $k$ is the {\em period}
of $\rho$. A Nielsen path $\rho$ is called {\em indivisible} if it cannot
be expressed as the concatenation of shorter Nielsen paths. A path $\rho$
is a {\em pre-Nielsen path} if $f^k_\#(\rho)$ is Nielsen for some $k\geq 0$.

A decomposition of a path $\rho=\rho_1\cdot\rho_2\ldots\cdot\rho_s$ into
subpaths is called a {\em $k$-splitting} if
$f^k_\#(\rho)=f^k_\#(\rho_1)\cdots f^k_\#(\rho_s).$ Such a decomposition
is a {\em splitting} if it is a $k$-splitting for all $k>0$. We will
also use the notion of $k$-splittings of circuits
$\sigma=\rho_1\cdot\rho_2\ldots\cdot\rho_s$, which requires, in
addition, that there be no cancellation between $f^k_\#(\rho_s)$ and
$f^k_\#(\rho_1)$.

The following theorem was proved in \cite{hb1}.
\begin{thm}[{\cite[Theorem 5.12]{hb1}}]\label{rtt}
Every outer automorphism $\mathcal O$ of $F$ is represented by a
homotopy equivalence $f\co G\rightarrow G$ such that each exponentially
growing stratum $H_r$ has the following properties:
\begin{enumerate}
\item If $E$ is an edge in $H_r$, then the
first and last edges in $f(E)$ are contained in $H_r$.
\item If $\beta$ is a nontrivial path in $G_{r-1}$ with endpoints in
$G_{r-1}\cap H_r$, then $f_\#(\beta)$ is nontrivial.\label{lowerstrat}
\item If $\rho$ is an $r$-legal path, then $f_\#(\rho)$ is an $r$-legal path.
\end{enumerate}
\end{thm}
We call $f$ a {\em relative train track map}.

%An outer automorphism $\mathcal O$ of $F$ is called {\em reducible} if it
%preserves the conjugacy class of a proper free factor of $F$. $\mathcal O$
%is called {\em irreducible} if it is not reducible. If $\mathcal O$ is
%irreducible, then it has a relative train track representative
%$f\co G\rightarrow G$ whose filtration has only one nonempty element $H_1=G$,
%with irreducible transition matrix. The properties of relative train tracks
%show that for every edge $E$ of $G$, the image $f^n(E)$ is an immersion
%for all $n>0$. In this case, we call $f$ a {\em train track map} (or 
%{\em absolute} train track map), and we denote the growth rate of
%$H_1=G$ by $\lambda$.

A path $\rho$ in $G$ is said to be of {\em height $r$} if $\rho\subset G_r$
and $\rho\not\subset G_{r-1}$. If $H_r=\{E_r\}$ is a polynomially
growing stratum, then {\em basic paths} of height $r$ are of the form
$E_r\gamma$ or $E_r\gamma E_r^{-1}$, where $\gamma$
is a path in $G_{r-1}$. If $\tau$ is a closed Nielsen path in $G_{r-1}$
and $f(E_r)=E_r\tau^l$ for some $l\in\mathbb Z$, then paths of the form
$E_r\tau^k$ and $E_r\tau^kE_r^{-1}$ are {\em exceptional paths} of
height $r$.  Moreover, if $s<r$, $\tau\subset G_{s-1}$, and
$f(E_s)=E_s\tau^m$, then $E_r\tau^kE_s^{-1}$ is also a exceptional path
of height $r$.

For our purposes, the properties of relative train track maps are not strong
enough, so we will use the notion of improved train track maps
constructed in \cite{tits1}. We only list the properties used in this
paper.

\begin{thm}[{\cite[Theorem~5.1.5, Lemma~5.1.7, and Proposition~5.4.3]{tits1}}]
\label{ttimproved}
For every outer automorphism $\mathcal O$ of $F$, there exists
an exponent $k>0$ such that $\mathcal O^k$ is represented by a relative
train track map $f\co G\rightarrow G$ with the following additional properties:
\begin{enumerate}
\item \label{zerostrat} If $H_{r}$ is a zero stratum, then $H_{r+1}$ is an
exponentially growing stratum, and the restriction of $f$ to $H_r$ is an
immersion. $H_r$ is a zero stratum if and only if it is
the union of the contractible components of $G_r$.
\item\label{fixedvert} If $v$ is a vertex, then $f(v)$ is a fixed
vertex. If $H_r$ is a polynomially growing stratum and $G'$ is the
collection of noncontractible components of $G_{r-1}$, then all vertices
in $H_r\cap G'$ are fixed.
\item\label{npprop} If $H_r$ is an exponentially growing stratum, then
there is at most one indivisible Nielsen path $\tau$ of height $r$.  If
$\tau$ is not closed and if it starts and ends at vertices, then at
least one endpoint of $\tau$ is not contained in $H_r\cap G_{r-1}$.
\item \label{polyim} If $H_r$ is a polynomially growing stratum,
then $H_r$ consists of
a single edge $E_r$, and $f(E_r)=E_r\cdot u_r$ for some closed path
$u_r\subset G_{r-1}$ whose base point is fixed by $f$.

If $\sigma\subset G_r$ is a basic path of height $r$
that does not split as a concatenation of two basic paths of height $r$
or as a concatenation of a basic path of height $r$ with a path contained
in $G_{r-1}$, then either $f^k_\#(\sigma)=E_r\cdot\sigma'$ for some
$k\geq 0$, or $u_r$ is a Nielsen path and $f^k_\#(\sigma)$ is an
exceptional path of height $r$ for some $k\geq 0$.

\end{enumerate}
\end{thm}
We call $f$ an {\em improved} relative train track map.

Finally, we state a lemma from \cite{tits1} that simplifies the study of
paths intersecting strata of polynomial growth.

\begin{lem}[{\cite[Lemma 4.1.4]{tits1}}]\label{pgsplit}
Let $f\co G\rightarrow G$ be an improved train track map with a
polynomially growing stratum $H_r$. If $\rho$ is a path in $G_r$, then
it splits as a concatenation of basic paths of height $r$ and paths in
$G_{r-1}$.
\end{lem}

\begin{rem}\label{howtosplit}
In fact, part \ref{polyim} of Theorem \ref{ttimproved} implies that
subdividing $\rho$ at the initial endpoints of all occurrences of $E_r$
and at the terminal endpoints of all occurrences of $E_r^{-1}$ yields a
splitting of $\rho$ into basic paths of height $r$ and paths in
$G_{r-1}$.
\end{rem}

Observe that if $H_r=\{E_r\}$ is a polynomially growing stratum, then
$f^k_\#(E_r)=E_r\cdot u_r \cdot f_\#(u_r) \cdot \ldots \cdot f^{k-1}_\#(u_r)$.
Each subpath of the form $f^i_\#(u_r)$ is called a {\em block} of
$f^k_\#(E_r)$. Since there is no cancellation between successive blocks, it
makes sense to refer to the infinite path
\begin{equation}
	R_r=u_r \cdot f_\#(u_r) \cdot f_\#^2(u_r) \cdot \ldots \label{raydef}
\end{equation}
as the {\em eigenray} of $E_r$.

\begin{rem}[A note on terminology]\label{terminology}
The notion of a {\em polynomially growing stratum} $H_r=\{E_r\}$ first
appeared in \cite{hb1}. Polynomially growing strata are called {\em
nonexponentially growing} strata in \cite{tits1}. Both terms are
somewhat misleading because the function $k\mapsto |f_\#^k(E_r)|$ may
grow exponentially (see \lemref{fastpoly}).
\end{rem}

Given an improved train track map $f\co G\rightarrow G$, we construct a metric
on $G$. If $H_r$ is an exponentially growing stratum, then its transition
matrix $M_r$ has a unique positive left eigenvector $v_r$
(corresponding to $\lambda_r$) whose smallest entry equals one
\cite{gantmacher}.
For an edge $E_i$ in $H_r$, the eigenvector $v_r$ has an entry $l_i>0$
corresponding to $E_i$. We choose a metric on $G$ such that $E_i$ is
isometric to an interval of length $l_i$, and such that edges in zero strata
or in polynomially growing strata are isometric to an interval of length
one. For a path $\rho$, we denote its length by $\mathcal{L}(\rho).$
Note that if the endpoints of $\rho$ are vertices, then the number of
edges in $\rho$ provides a lower bound for $\mathcal{L}(\rho)$.
Moreover, if $f$ is an absolute train track map, then $f$ expands the
length of legal paths by the factor $\lambda$.

\begin{rem}
We merely choose this metric for convenience. All statements here are
invariant under bi-Lipschitz maps, but our metric of choice simplifies
the presentation of our arguments.
\end{rem}

We are now ready to state the main technical result of this paper.

\begin{thm}\label{mainthm}
Let $\phi\co F\rightarrow F$ be an an automorphism. Then there exists
an improved relative train track map representing some positive power of 
$\phi$ for which there exists a constant $K\geq 1$ with the following
property: For any pair of exponents $N,i$ satisfying $0\leq i \leq N$,
the following two statements hold:
\begin{enumerate}
\item If $\sigma$ is a circuit in $G$, then
\[
    \mathcal{L}\left(f^i_\#(\sigma)\right)\leq
    K\left(\mathcal{L}(\sigma)+\mathcal{L}\left(f^N_\#(\sigma)\right)\right).
\]
\item If $\rho$ is a path in $G$ that starts and ends at vertices, then
\[
    \mathcal{L}\left(f^i_\#(\rho)\right) \leq
    K\left(\mathcal{L}(\rho)+\mathcal{L}\left(f^N_\#(\rho)\right)\right).
\]
\end{enumerate}
Given the improved relative train track map $f\co G\rightarrow G,$ the
constant $K$ can be computed.
\end{thm}

We will present the proof of \thmref{mainthm} in \secref{polysec} and
\secref{pfsec}.  Right now, we show how \thmref{mainapp} follows from
\thmref{mainthm}.

\begin{proof}[Proof of \thmref{mainapp}]
Let $\phi\co F\rightarrow F$ be an automorphism of a finitely generated
free group $F=\langle x_1,\ldots,x_n\rangle$.
The first part of \thmref{mainthm} immediately implies that
the first part of \thmref{mainapp} holds for some positive power $\phi^k$,
i.e., there exists some $K'\geq 1$ such that for all $0\leq i\leq N$ and
$w\in F$, we have
\[
    ||\phi^{ik}_\#(w)||\leq K'\left(||w||+||\phi^{Nk}_\#(w)||\right),
\]
where we compute lengths with respect to the generators $x_1,\ldots,x_n$.

Let $L=\max\{|\phi(x_i)|\}$. Then, for $0\leq j<k$, we have
\[
	L^{-k}||\phi^{ik+j}(w)||\leq ||\phi^{ik}(w)||\leq L^k||\phi^{ik+j}(w)||
\]
for all $w\in F$. We conclude that for all $0\leq i\leq N$ and $w\in F$,
we have
\[
    L^{-k}||\phi^{i}_\#(w)||\leq K'\left(||w||+L^k||\phi^N_\#(w)||\right),
\]
so that the first part of \thmref{mainapp} holds with $K=L^{2k}K'$.

In order to prove the second assertion, we modify a trick from \cite{bfgafa2}.
Let $F'$ be the free group generated by $x_1,\ldots,x_n$ and an additional
generator $a$. We define an automorphism $\psi\co F'\rightarrow F'$ by
letting $\psi(x_i)=\phi(x_i)$ for all $1\leq i\leq n$, and $\psi(a)=a$.

By the previous step, the first part of \thmref{mainapp} holds for $\psi$,
with some constant $K'\geq 1$. Let $w$ be some word in $F$. Then, for all
$i\geq 0$, $\psi^i(aw)$ is a cyclically reduced word in $F'$, so that we
have $|\phi^i(w)|+1=||\psi^i(aw)||$. We conclude that
\[
	|\phi^i(w)|+1\leq K'(|w|+|\phi^N(w)|+2),
\]
for all $0\leq i\leq N$. Now the second assertion of \thmref{mainapp}
holds with $K=2K'$.
\end{proof}

\section{More on train tracks}\label{train2}

Thurston's bounded cancellation lemma is one of the fundamental tools in this
paper. We state it in terms of homotopy equivalences of graphs.

\begin{lem}[Bounded cancellation lemma \cite{cooper}]\label{bcc}
Let $f:G\rightarrow G$ be a homotopy equivalence. There exists a constant
$\mathcal{C}_f$, depending only on $f$, with the property that for any
tight path $\rho$ in $G$ obtained by concatenating two paths $\alpha,
\beta$, we have
\[
\mathcal{L}(f_\#(\rho))\geq \mathcal{L}(f_\#(\alpha)) +
\mathcal{L}(f_\#(\beta)) - \mathcal{C}_f.
\]
\end{lem}

An upper bound for $\mathcal{C}_f$ can easily be read off from the map $f$
\cite{cooper}.  Let $f:G\rightarrow G$ be an improved relative train
track map with an exponentially growing stratum $H_r$ with growth rate
$\lambda_r$.  The {\em $r$-length} of a path $\rho$ in $G$,
$\mathcal{L}_r(\rho)$, is the total length of $\rho\cap H_r$.

If $\beta$ is an $r$-legal path in $G$ whose $r$-length satisfies
$\lambda_r \mathcal{L}_r(\beta)-2\mathcal{C}_f>\mathcal{L}_r(\beta)$ and
$\alpha, \gamma$ are paths such that the concatenation
$\alpha\beta\gamma$ is an immersion, then the $r$-length of the segment
in $f^k_\#(\alpha\beta\gamma)$ corresponding to $\beta$
(\convref{howtotighten}) will tend to infinity as $k$ tends to infinity.
The {\em critical length} $\mathcal{C}_r$ of $H_r$ is the infimum of the
lengths satisfying the above inequality, i.\ e.,
\begin{equation}
	\mathcal{C}_r=\frac{2\mathcal{C}_f}{\lambda_r-1}. \label{critlen}
\end{equation}

We now list some additional technical results about improved train
track maps. The following lemma is an immediate consequence of
\cite[Proposition~6.2]{pbhyp}. If $H_r$ is an exponentially growing
stratum, and $\rho$ is a path of height $r$, we let $n(\rho)$ denote
the number of $r$-legal segments in $\rho$.

\begin{lem}\label{tricho}
Let $f:G\rightarrow G$ be a relative train track map, and
let $H_r$ be an exponentially growing stratum. For each $L>0$, there
exists some computable exponent $M>0$ such that if $\rho$ is a path or
circuit in $G_r$ containing at least one full edge in $H_r$, one of
the following three statements holds:
\begin{enumerate}
\item $f^M_\#(\rho)$ has an $r$-legal segment of $r$-length greater than $L$.
\item $n(f^M_\#(\rho))<n(\rho)$.
\item $\rho$ can be expressed as a concatenation $\tau_1\rho'\tau_2$,
where $\tau_1, \tau_2$ each contain at most one $r$-illegal turn,
the $r$-length of the $r$-legal segments of $\tau_1, \tau_2$ is at most $L$,
and $\rho'$ splits as a concatenation of pre-Nielsen paths
(with one $r$-illegal turn each) and segments in $G_{r-1}$. \label{case3}
Moreover, $f^M_\#(\rho')$ is a concatenation of Nielsen paths of height
$r$ and segments in $G_{r-1}$.
\end{enumerate}
\end{lem}

\begin{rem}\label{tricho2}
\
\begin{itemize}
\item The statement of \lemref{tricho} in \cite{pbhyp} does not
explicitly mention the computability of $M$. The proof, however,
only uses counting arguments, from which the constant $M$ can be
computed.
\item The presence of the subpaths $\tau_1, \tau_2$ in \partref{case3} is
an artifact of the fact that $\rho$ need not start or end at fixed points
if it is a path. If $\rho$ starts at a fixed point, then $\tau_1$ will be
trivial, and if $\rho$ ends at a fixed point, then $\tau_2$ will be trivial.
\item The actual statement of \cite[Proposition~6.2]{pbhyp} does not mention
circuits since they were not a concern in the context of \cite{pbhyp}. The
proof, however, works for circuits as well as paths. If the first two
statements of \lemref{tricho} do not hold, than the third statement will
hold with $\tau_1$ and $\tau_2$ trivial.
\end{itemize}
\end{rem}

From now on, we assume that $f\co G\rightarrow G$ that $f$ is an
improved train track map. Throughout the rest of this section, let $M$
be the constant from \lemref{tricho} for some fixed $L>\mathcal{C}_r$
(\eqref{critlen}).

Let $H_r=\{E_r\}$ be a polynomially growing stratum. We say that $H_r$
is {\em truly polynomial} if $u_r$ is trivial or, inductively, if $u_r$
is a concatenation of truly polynomial edges and Nielsen paths in
exponentially growing strata. Clearly, if $E_r$ is truly polynomial,
then the map $k\mapsto |f_\#^k(E_r)|$ grows polynomially.  We say that a
polynomially growing stratum is {\em fast} if it is not truly
polynomial.

The following lemma give us an understanding of the growth of fast
polynomial strata.

\begin{lem}\label{fastpoly}
There exists an exponent $k_0$ with the following property: For all fast
polynomial strata $H_r=\{E_r\}$ there exists some $s<r$ such that $H_s$
is of exponential growth and $f_\#^{k_0}(E_r)$ contains an $s$-legal
subpath of height $s$ whose $s$-length exceeds $\mathcal{C}_s$.
\end{lem}

In particular, this lemma implies that fast polynomial strata grow
exponentially. Given an improved relative train track map, we can find
$k_0$ by successively evaluating $f_\#, f_\#^2, \ldots$ until we see
long legal segments in all images of fast polynomial edges.

\begin{proof}
We introduce {\em classes} of fast polynomial edges. Let $H_r=\{E_r\}$
be a fast polynomial edge such that $f(E_r)=E_ru_r$.  We say that $H_r$
has {\em class~$1$} if there exists some $s<r$ such that $H_s$ is an
exponentially growing stratum, $u_r\cap H_s$ does not only consist of
Nielsen paths and paths of height less than $s$, and if $u_r$ contains
any polynomial edges $E_t$ for some $t>s$, then $E_t$ is truly polynomial.
We recursively define a fast polynomial edge $E_r$ to have {\em class~$k$}
if the highest class of edges in $u_r$ is $k-1$.

If $H_r$ has class~$1$, then $u_r$ contains a subpath $\rho$ of height
$s$ such that $f_\#^k(\rho)$ contains a long $s$-legal segment for some
sufficiently large $k$ (\lemref{tricho}).  If $u_r$ contains any
subpaths whose height exceeds $s$, then by definition those subpaths
will grow at most polynomially, so that eventually, the exponential
growth of $\rho$ will prevail.

In order to prove the lemma for an edge of class~$k$, $k>1$, we observe
that no edges of class~$k-1$ are cancelled when $f^m(u_r)$ is tightened
to $f_\#^m(u_r)$. Now the lemma follows by \thmref{ttimproved},
\partref{polyim}, and induction.
\end{proof}

Assume that $H_r$ is an exponentially growing stratum, and let $\rho$
be a path of height $r$. If $H_r$ does not support a closed Nielsen
path, then we let $N(\rho)=n(\rho)$. If $H_r$ supports a closed
Nielsen path, then we let $N(\rho)$ equal the number of legal segments
in $\rho$ that do not overlap with a Nielsen subpath of $\rho$.

The following lemma is a generalization of \cite[Lemma~6.4]{pbhyp}.
\begin{lem}\label{expdecay}
Assume that $H_r$ is an exponentially growing stratum.
There exist computable constants $\lambda>1, N_0$ with the following
property: If $f^M_\#(\rho)$ does not contain a legal segment of length
at least $L$, and if $N(\rho)>N_0$, then
\[
	N(f^M_\#(\rho))\leq \lambda^{-1}N(\rho).
\]
Regardless of $N(\rho)$, we have
\[
	N(f^M_\#(\rho))\leq \lambda^{-1}N(\rho)+1.
\]
\end{lem}

\begin{proof}
If $H_r$ does not support a closed Nielsen path, then the proof of
\cite[Lemma~6.4]{pbhyp} goes through unchanged. We repeat the argument
here because the ideas of the proof show up more clearly in this case.

If $H_r$ does not support a closed Nielsen path, then the proof is
based on the following observation: If $N(\rho)=6$
and $f_\#(\rho)$ does not contain a long legal segment, then
$N(f^M_\#(\rho))<6$. Suppose otherwise, i.e., $N(f^M_\#(\rho))=N(\rho)$.
Then, by \lemref{tricho}, $f^M_\#(\rho)=\tau_1\gamma\tau_2$, where
$\gamma$ is a concatenation of three indivisible Nielsen paths of
height $r$ and paths in $G_{r-1}$. This is impossible because by
\thmref{ttimproved}, \partref{npprop}, we can concatenate no more
than two indivisible Nielsen paths of height $r$ with paths in $G_{r-1}$.

Hence, of every six consecutive legal segments in $\rho$, at least
one cancels completely when $f^M(\rho)$ is tightened to $f^M_\#(\rho)$.
This implies that if $N(\rho)\geq 6$, then
$N(f^M_\#(\rho))\leq \frac{10}{11}N(\rho)$. In order to see
why this choice of $\lambda$ works, we just observe that if $\rho$ consists
of eleven legal segments and the sixth one cancels in $f^M_\#(\rho)$,
then there are no six consecutive legal segments that survive in
$f^M_\#(\rho)$.

This completes the proof of the first inequality, with $\lambda=\frac{11}{10}$
and $N_0=5$, if $H_r$ does not support a closed Nielsen path. Regarding the
second inequality, we remark that if $N(\rho)\leq N_0$, then
$N(f^M_\#(\rho))\leq N(\rho)\leq \lambda^{-1}N(\rho)+1$.

We now assume that $H_r$ supports a closed indivisible Nielsen path
$\sigma$.  The proof in this case is based on the following consequence
of \lemref{tricho}. If $\gamma$ a path of height $r$,
$n(\gamma)=n(f^M_\#(\gamma))=4$, and $f^M_\#(\gamma)$ does not contain
a long legal segment, then
$f^M_\#(\gamma)=\tau_1\sigma^{\pm 1}\tau_2$, where $\tau_1$ and $\tau_2$
are as in \lemref{tricho}. Intuitively, this means that if few legal
segments disappear, then many Nielsen paths will appear. Since $N(\rho)$
only counts those legal segments that do not overlap with a Nielsen path,
this observation will yield the desired estimate.

First, consider a path $\gamma$ of height $r$ that does not contain any
Nielsen subpaths, i.e., we have $N(\gamma)=n(\gamma)$. If $N(\gamma)\geq 4$,
then for every four consecutive legal segments whose images do not cancel
completely in $f^M_\#(\gamma)$, $f^M_\#(\gamma)$ contains at least one
Nielsen subpath, so that we have $N(f^M_\#(\gamma))\leq \frac{6}{7}N(\gamma)$,
using the same reasoning as above.

We claim that if $\gamma$ starts and ends at fixed points, then, by
\remref{tricho2}, we have $N(f^M_\#(\gamma))\leq \frac{6}{7}N(\gamma)$
regardless of $N(\gamma)$. To this end, we first argue that if $\gamma$
starts and ends at fixed points, then $n(f^M_\#(\gamma))<n(\gamma)$. If
this were not true, then, by \lemref{tricho} we would have
$f^M_\#(\gamma)=\sigma^m$ for some $m\in\mathbb{Z}$, which would imply
that $\gamma=\sigma^m$ because $\gamma$ starts and ends at fixed points.
This is a contradiction since we assumed that $\gamma$ does not contain
any Nielsen subpaths. Now, if $n(\gamma)=N(\gamma)<4$, then we conclude that
$N(f^M_\#(\gamma))\leq n(f^M_\#(\gamma))<n(\gamma)$. Now $n(\gamma)<4$
implies that $\frac{6}{7}n(\gamma)\geq n(\gamma)-1$, which implies
that $N(f^M_\#(\gamma))\leq \frac{6}{7}n(\gamma)=\frac{6}{7}N(\gamma)$.

After these preparations, we express $\rho$ as a concatenation
\[
	\rho=\rho_1\sigma^{n_1}\rho_2\sigma^{n_2}\rho_3\cdots\rho_{k}\sigma^{n_k}
		\rho_{k+1},
\]
where $n_1,\ldots,n_k\in\mathbb{Z}$, and none of the subpaths $\rho_i$
contains a Nielsen subpath.

Note that the subpaths $\rho_2,\ldots,\rho_k$ start and end at the
base point $v$ of the Nielsen path $\sigma$, which is fixed by $f$.
Hence, for $2\leq i\leq k$, we have
$N(f^M_\#(\rho_i))\leq \frac{6}{7}N(\rho_i)$, and we have
$N(f^M_\#(\rho_1))\leq \frac{6}{7}N(\rho_1)$ (resp.\ 
$N(f^M_\#(\rho_{k+1}))\leq \frac{6}{7}N(\rho_{k+1})$) if
$N(\rho_1)\geq 4$ (resp.\ $N(\rho_{k+1})\geq 4$).

If $N(\rho_1)<4$ and $N(\rho_{k+1})<4$, we have
\begin{eqnarray*}
	N(f^M_\#(\rho)) & \leq &
		N(\rho_1)+N(\rho_{k+1})+\frac{6}{7}(N(\rho)-N(\rho_1)-N(\rho_{k+1})) \\
		& \leq & 6+\frac{6}{7}(N(\rho)-6)
		\leq \frac{6}{7}(1+N(\rho)).
\end{eqnarray*}
Similar estimates yield that $N(f^M_\#(\rho))\leq \frac{6}{7}(1+N(\rho))$
regardless of $N(\rho_1)$ and $N(\rho_{k+1})$.

If $N(\rho)>11$, then $\frac{6}{7}(1+N(\rho))\leq \frac{13}{14}N(\rho)$,
which implies that $N(f^M_\#(\rho))\leq \frac{13}{14}N(\rho)$ if
$N(\rho)>11$, so that the first inequality of the lemma holds with
$\lambda=\frac{14}{13}$ and $N_0=11$. As for the second inequality,
we remark that $N(f^M_\#(\rho))\leq N(\rho)$ and, if $N(\rho)\leq 11$,
then $N(\rho)\leq \lambda^{-1}N(\rho)+1$.
\end{proof}

The next lemma is a statement about the (absence of) cancellation between
eigenrays of polynomially growing strata. It is a stronger version of
\cite[Sublemma~1, Page~587]{tits1}.

\begin{lem}\label{sublemma}
Let $H_i=\{E_i\}$ and $H_j=\{E_j\}$ be polynomially growing strata.  Let
$S_i$ (resp.\ $S_j$) be an initial segment of $E_iR_i$ (resp.  $E_jR_j$,
see \eqref{raydef}) such that the concatenation $S_i\bar{S}_j$ is a
path.  If $E_i$ grows faster than linearly and if an entire block of
$R_j$ is canceled in $f^k_\#(S_i\bar{S}_j)$ for some $k\geq 0$, then no
entire block of $R_i$ will be canceled in $f^l_\#(S_i\bar{S}_j)$ for any
$l\geq 0$.
\end{lem}

\begin{proof}
Suppose that at least one block of both $S_i$ and $S_j$ cancels.
Then there are paths $\alpha$, $\beta$, and $\gamma$ such that
$f^k_\#(u_i)=\beta\gamma$, $f^l_\#(u_j)=\alpha\beta$ for some
$k, l\geq 0$, and $f_\#(\alpha)=\gamma$ (see \figref{slsplitfig}).

\begin{figure}
\renewcommand{\epsfsize}[2]{0.5\textwidth}
\begin{center}
\psfrag{Ei}{$E_i$}
\psfrag{Ej}{$E_j$}
\psfrag{a}{$\alpha$}
\psfrag{b}{$\beta$}
\psfrag{c}{$\gamma$}
\psfrag{fkui}{$f^k_\#(u_i)$}
\psfrag{fluj}{$f^l_\#(u_j)$}
\epsfbox{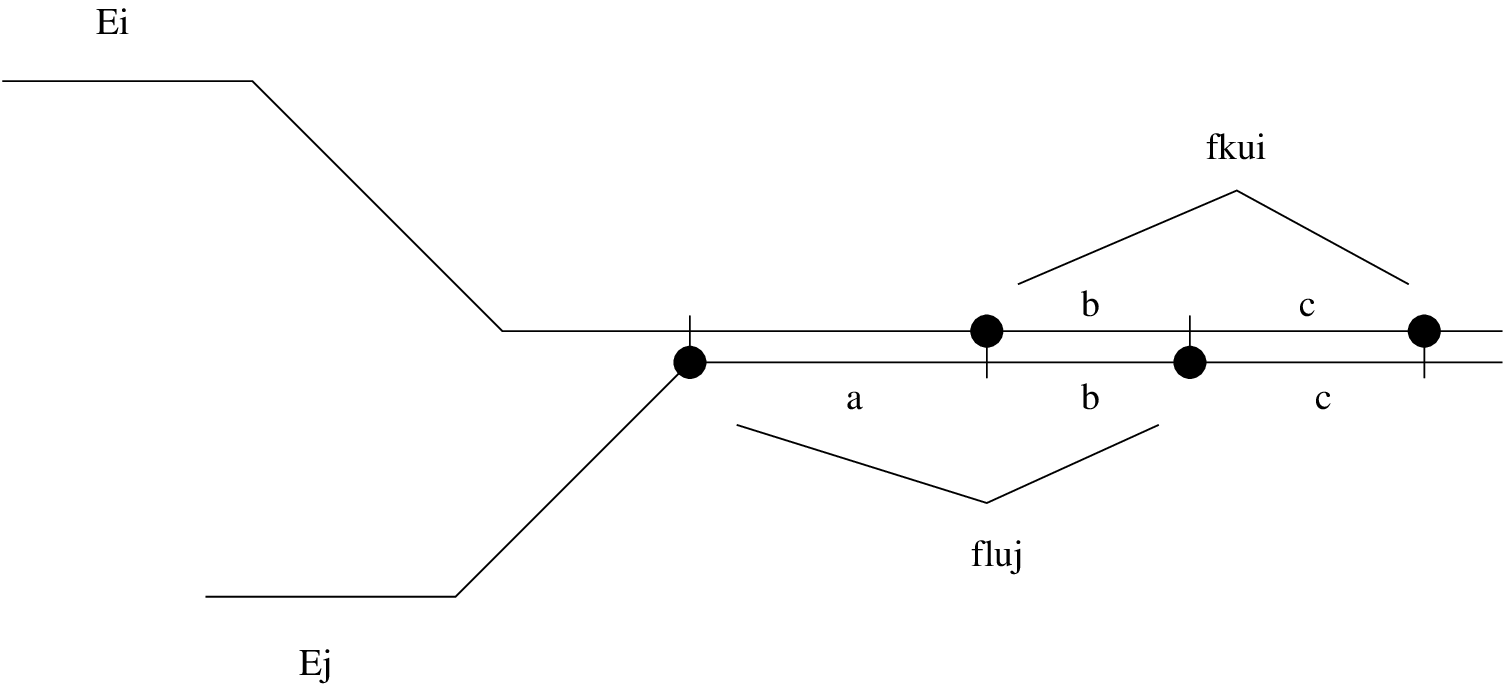}
\end{center}
\caption{The idea of the proof of \lemref{sublemma}.}
\label{slsplitfig}
\end{figure}

In particular, we have
\[
    R_i=u_if_\#(u_i)\cdots f_\#^{k-1}(u_i)\beta f_\#(\alpha) f_\#(\beta)
    f^2_\#(\alpha) \ldots,
\]
and
\[
	R_j=u_jf_\#(u_j)\cdots f_\#^{l-1}(u_j)\alpha \beta f_\#(\alpha) f_\#(\beta) f^2_\#(\alpha) \ldots.
\]
In particular, the path
$\rho=E_iR_i^{k-1}\bar{\alpha}\bar{R}_j^{l-1}\bar{E}_j$ does not split.
By \thmref{ttimproved}, $\rho$ is a exceptional path, and both $E_i$ and
$E_j$ grow linearly.
\end{proof}

\section{Terminology and examples}\label{exsec}

In this section, we discuss some examples that illustrate some of the
main issues that we need to address in the proof of \thmref{mainthm}.
Although we are not primarily concerned with free-by-cyclic groups
in this article, the language of free-by-cyclic groups will streamline
the exposition.

Given a free group $F_n=\langle x_1,\ldots,x_n \rangle$ and an
automorphism $\phi$ of $F_n$, the {\em mapping torus} of $\phi$
is the free-by-cyclic group
\[
	M_\phi=\langle x_1,\ldots,x_n,t~|~t^{-1}x_it\phi(x_i^{-1})\rangle.
\]
The letter $t$ is called the {\em stable letter} of $M_\phi$.

A reduced word $w$ in the generators of $M_\phi$ is a {\em hallway}
if $w$ represents the trivial element of $M_\phi$ and if $w$ can
be expressed as $w=w_1w_2$ such that $w_1$ only contains negative
powers of $t$ and $w_2$ only contains positive powers of $t$ \cite{combi}.
Hallways of the form $t^{-k}xt^k\phi^k(x^{-1})$, for $x\in F_n$,
are said to be {\em smooth}.

\begin{figure}
\renewcommand{\epsfsize}[2]{0.5\textwidth}
\begin{center}
\psfrag{k}{$k$}
\psfrag{w0}{$w_0$}
\psfrag{wi}{$w_i$}
\psfrag{wkp1}{$w_k$}
\psfrag{u1}{$u_1$}
\psfrag{uk}{$u_{k-1}$}
\psfrag{v1}{$v_1$}
\psfrag{v2}{$v_2$}
\psfrag{vk}{$v_{k-1}$}
\epsfbox{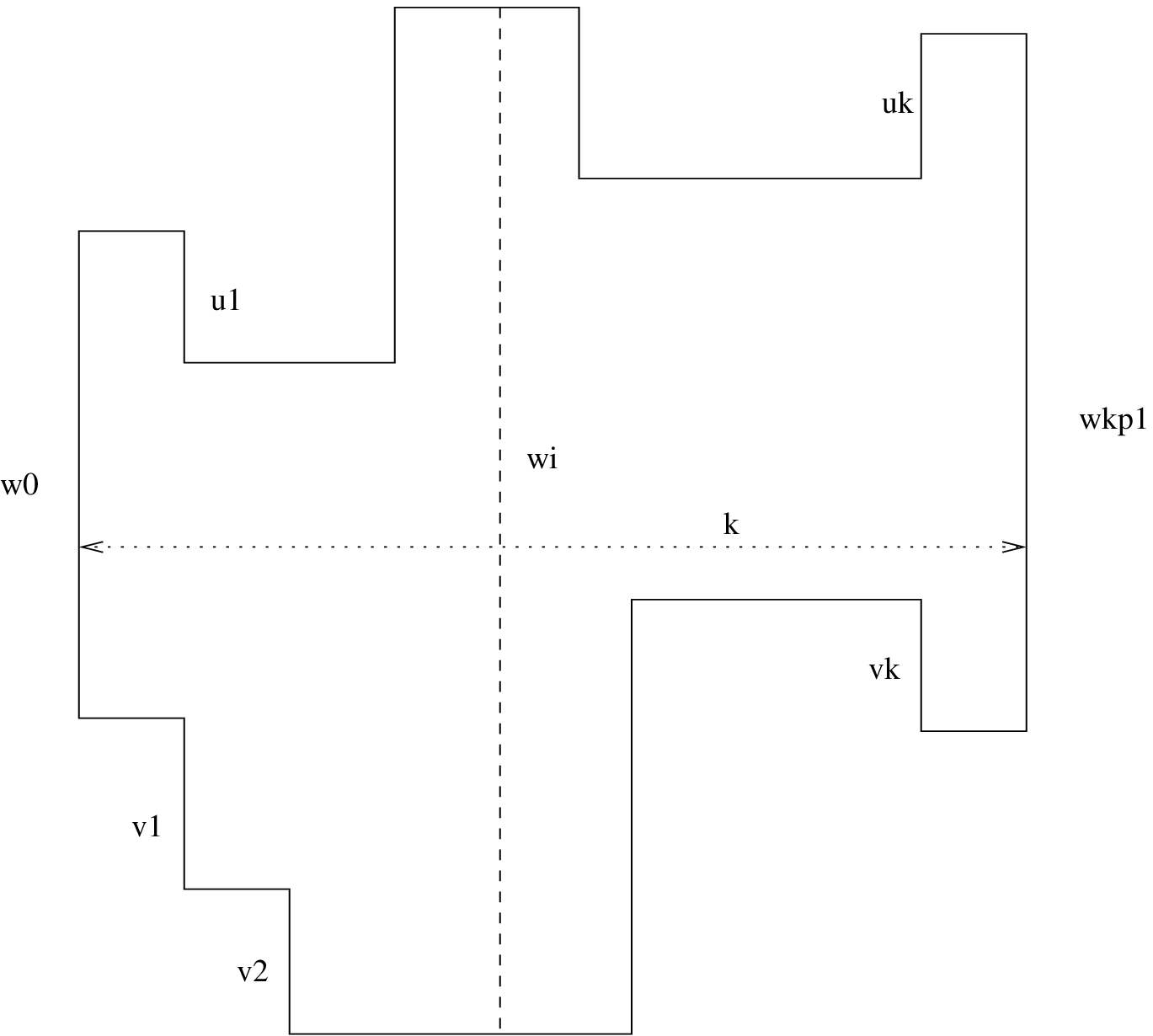}
\end{center}
\caption{A hallway.}
\label{hallwayfig}
\end{figure}

Any hallway $w$ can be expressed as
\[
	w=t^{-1}u_{k-1}t^{-1}u_{k-2}t^{-1}\cdots t^{-1}u_1t^{-1}w_0
		tv_1tv_2t\cdots tv_{k-1}tw_k^{-1},
\]
where $w_0,w_k,u_1,\ldots,u_{k-1},v_1,\ldots,v_{k-1}$ are elements
of $F_n$. The words $u_i$ and $v_i$ may be empty. In fact, a hallway
is smooth if and only if all the $u_i$ and $v_i$ are trivial.
For $1\leq i< k$, we define $w_i$ to be the word obtained by
tightening $u_i\phi(w_{i-1})v_i$.  Since $w$ represents the identity,
we have $w_k=\phi(w_{k-1})$.  We call $w_i$ the $i$-th {\em slice} of
$w$. The number $k$ is the {\em duration} $\mathcal{D}(w)$ of the hallway.
\figref{hallwayfig} illustrates these notions.

We say that the instances of letters of $F_n$ that occur in the spelling
of $w$ are {\em visible}. \thmref{mainapp} states that if $w$ is a smooth
hallway, then the length of each $w_i$ is bounded by a constant multiple
of the number of visible edges in $w$.

The following examples illustrate the main issues that arise in the proof.
For the remainder of this section, let $F_6=\langle a,b,c,d,x,y \rangle$,
and define $\phi$ by letting
\begin{eqnarray*}
	a &\mapsto & a \\
	b &\mapsto & ba \\
	c &\mapsto & caa \\
	d &\mapsto & dc \\
	x &\mapsto & y \\
	y &\mapsto & xcy. \\
\end{eqnarray*}
This automorphism admits the stratification $H_1=\{a\}$, $H_2=\{b\}$,
$H_3=\{c\}$, $H_4=\{d\}$, and $H_5=\{x,y\}$. The restriction of $f$
to the filtration element $G_3=H_1\cup H_2\cup H_3$ grows linearly,
the restriction to $G_4$ grows quadratically, and the stratum
$H_5$ is of exponential growth.

The first example illustrates the behavior of smooth hallways in
linearly growing filtration elements.

\begin{figure}
\renewcommand{\epsfsize}[2]{0.4\textwidth}
\begin{center}
\psfrag{b}{$b$}
\psfrag{c}{$c$}
\psfrag{a2}{$a^2$}
\psfrag{a3}{$a^3$}
\psfrag{ak}{$a^k$}
\epsfbox{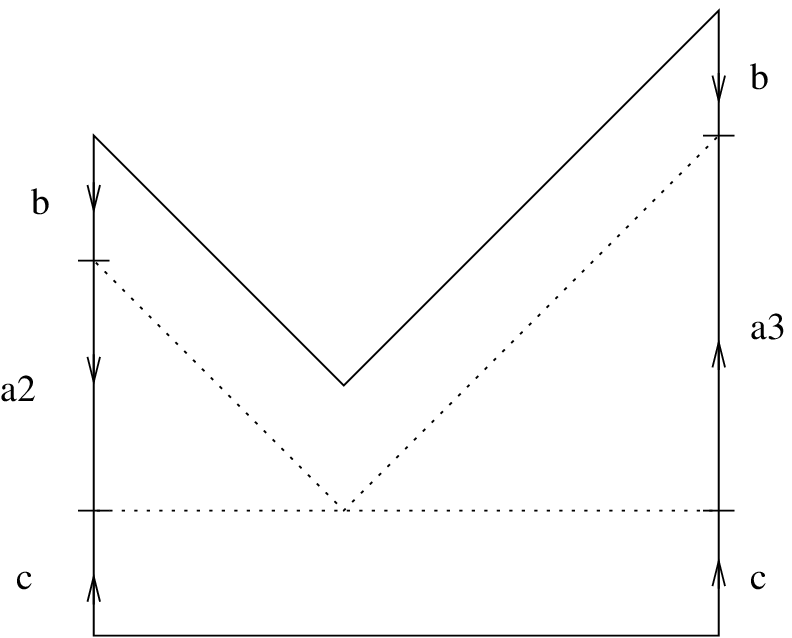}
\end{center}
\caption{Illustration of \exref{smoothex}.}
\label{smoothfig}
\end{figure}

\begin{ex}\label{smoothex}
Let $w_0$ be a word from the list $a^m,ba^mb^{-1},ca^mc^{-1}$, for some
integer $m$. Then $\phi(w_0)=w_0$, so that the length of any slice of
the hallway $t^{-k}w_0t^kw_0^{-1}$ is the same as the length of $w_0$.
Now, let $w_0$ be a word from the list $ba^m,ca^m,ca^mb^{-1}$.
If $m\geq 0$, then $|\phi^{k+1}(w_0)|=|\phi^k(w_0)|+1$ for any
$k\geq 0$. If $m<0$, then $|\phi^{k+1}(w_0)|=|\phi^k(w_0)|-1$ for
$0\leq k<-m$ (\figref{smoothfig}). Hence, the length of each slice of the
hallway $t^{-k}w_0t^k\phi^k(w_0^{-1})$ is bounded by the number of
visible letters.

If $w_0$ is an arbitrary word in $\langle a,b,c \rangle$,
then, by \remref{howtosplit},
it splits as a concatenation of words from the above lists and their
inverses, which implies that the lengths of slices of smooth hallways
is bounded by the number of visible letters, so that \thmref{mainapp}
holds with $K=1$.
\end{ex}

The next example shows that hallways that are not smooth may have slices
whose length is not bounded in terms of a constant multiple of the number
of visible edges.

\begin{figure}
\renewcommand{\epsfsize}[2]{0.8\textwidth}
\begin{center}
\psfrag{b}{$b$}
\psfrag{c}{$c$}
\psfrag{ak}{$a^k$}
\epsfbox{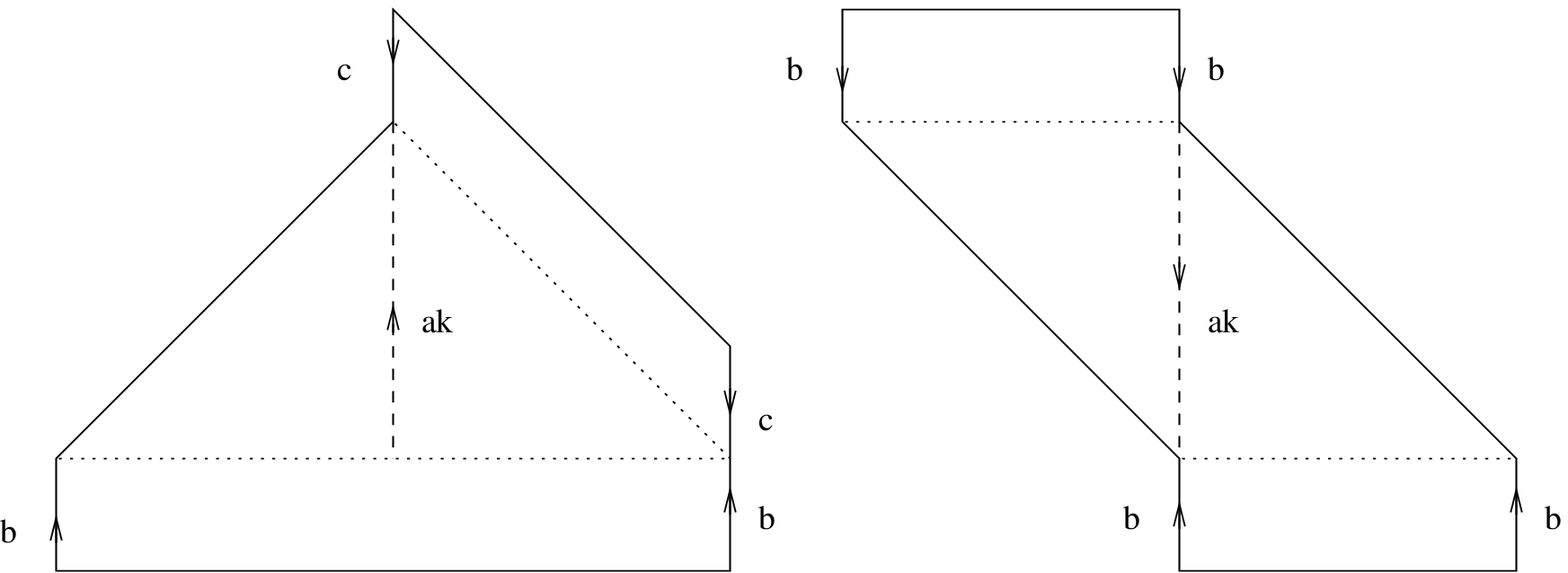}
\end{center}
\caption{Illustration of \exref{bulgeex}.}
\label{bulgefig}
\end{figure}

\begin{ex}\label{bulgeex}
Let $w=t^{-k}ct^{-k}b^{-1}t^{2k}bc^{-1}$. For $i<k$, we have
$w_i=a^{-i}b^{-1}$, and for $k\leq i \leq 2k$, we have
$w_i=ca^{2k-i}b^{-1}$ (\figref{bulgefig}). In particular, there is
a slice of length $k+2$ although there are only four visible edges in $w$.
Informally, one might say that hallways of this form bulge in the middle.
A similar bulge occurs for hallways of the form
$w=t^{-k}b^{-1}t^{-k}bt^kb^{-1}t^kb$.
\end{ex}

The next example shows that we need to control the size of such bulges
when proving \thmref{mainthm}.

\begin{figure}
\renewcommand{\epsfsize}[2]{0.4\textwidth}
\begin{center}
\psfrag{b}{$b$}
\psfrag{c}{$c$}
\psfrag{x}{$x$}
\psfrag{ak}{$a^k$}
\psfrag{pkxc}{$\phi^{-k}(xc)$}
\psfrag{pkx}{$\phi^k(x)$}
\epsfbox{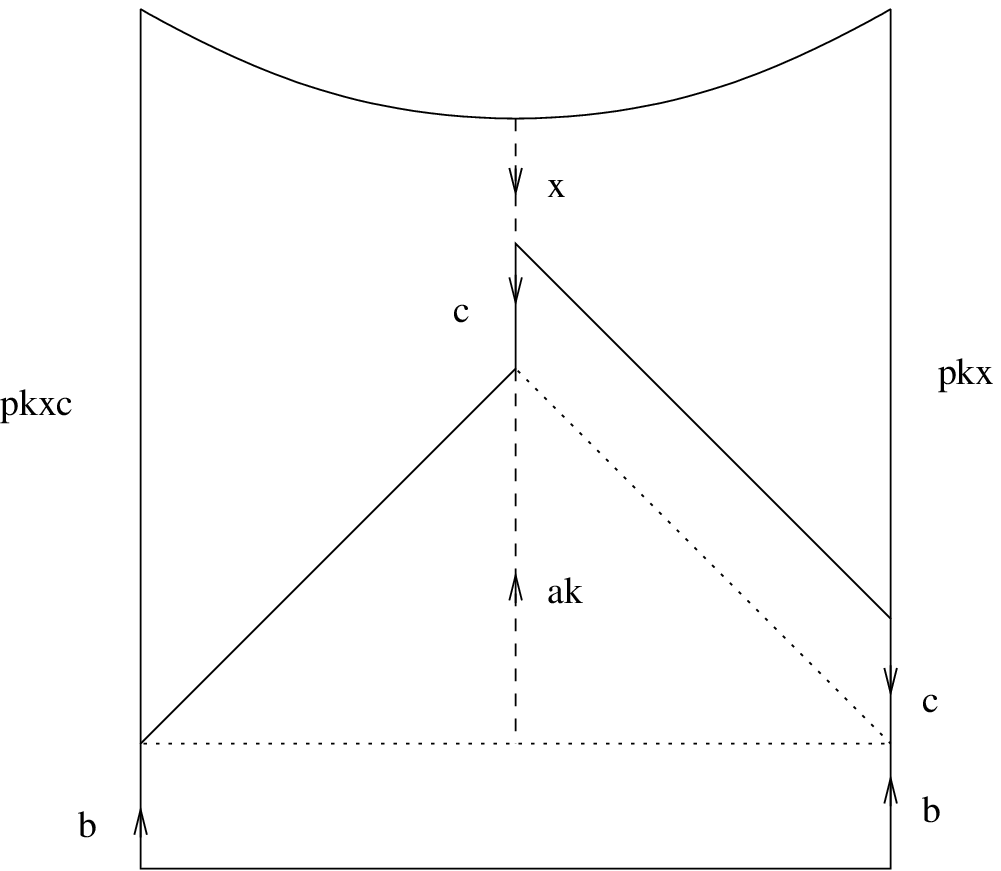}
\end{center}
\caption{Illustration of \exref{expex}.}
\label{expexfig}
\end{figure}

\begin{ex}\label{expex}
First, note that for $k\geq 1$, the last letter in the words
$f^{-k}(xc)$ is always one of $x,y,x^{-1},y^{-1}$, so that words of
the form $w_0=\phi^{-k}(xc)b^{-1}$ are reduced, and we have
$\phi^k(w_0)=xca^kb^{-1}$ and $\phi^{2k}(w_0)=\phi^k(x)b^{-1}$.

Hence, the smooth hallway $w=t^{-2k}w_0t^{2k}\phi^{2k}(w_0^{-1})$
contains a bulge like the first one in the previous example
(\figref{expexfig}). The presence of this bulge does not contradict
\thmref{mainapp} because $w$ contains a large number of visible
instances of the letters $x$ and $y$. This example shows that we cannot
consider the strata separately when proving \thmref{mainthm}.
\end{ex}

\begin{ex}\label{polyex}
If we let $w_0=\phi^{-k}(dc)b^{-1}$, then the smooth hallway
$w=t^{-2k}w_0t^{2k}\phi^{2k}(w_0^{-1})$ contains a bulge like
in \exref{expex}. This does not contradict \thmref{mainapp} as
$w$ contains a large number of visible instances of the letter $c$.
\end{ex}

Our final example illustrates a subtlety regarding linearly growing
strata.

\begin{ex}\label{linearsubtlety}
Let $F_4=\langle a,b,c,d \rangle$ and define $\psi$ by letting
\begin{eqnarray*}
	a & \mapsto & a \\
	b & \mapsto & ba \\
	c & \mapsto & ca \\
	d & \mapsto & dcb^{-1}.\\
\end{eqnarray*}
The map $\psi$ is a linearly growing automorphism, so in particular
the letter $d$ is of linear growth, although the image of $d$ contains
letters of linear growth other than $d$ itself.

Letters of linear growth may thus behave in two different ways; they
may contribute to the growth of images under successive applications
of $\psi$, or they may remain inert as parts of a fixed word. In the
proof of \thmref{mainthm}, we will need to distinguish letters of
linear growth according to their role.
\end{ex}

\begin{ex}\label{eglinear}
Let $F_3=\langle a, x, y \rangle$ and define $\phi$ by letting
\begin{eqnarray*}
	a & \mapsto & axyx^{-1}y^{-1} \\
	x & \mapsto & y^{-1} \\
	y & \mapsto & yx. \\
\end{eqnarray*}
The stratum $\{x, y\}$ grows exponentially, and we have
$\psi(xyx^{-1}y^{-1})=xyx^{-1}y^{-1}.$ This means that $a$ grows
linearly although it maps across an exponentially growing stratum.
This is another phenomenon that we need to consider when analyzing
strata of linear growth.
\end{ex}

The notion of hallways naturally extends to mapping tori of homotopy
equivalences of finite graphs. Specifically, a hallway $\rho$ in the
mapping torus of $f\co G\rightarrow G$ is a sequence of paths of the form
\[
	\rho=(\mu_{k-1}, \mu_{k-2}, \cdots , \mu_1, \rho_0
		, \nu_1, \nu_2, \cdots, \nu_{k-1}, \rho_k),
\]
where $\rho_0,\rho_k,\mu_1,\ldots,\mu_{k-1},\nu_1,\ldots,\nu_{k-1}$ are
paths in $G$, satisfying $f(\tau(\rho_0))=\iota(\nu_1),
f(\tau(\nu_i))=\iota(\nu_{i+1}), f(\tau(\nu_{k-1}))=\tau(\rho_k),
f(\iota(\rho_0))=\tau(\mu_1), f(\iota(\mu_i))=\tau(\mu_{i+1}),$ and
$f(\iota(\mu_{k-1}))=\iota(\rho_k),$ where $\iota(.)$ is the initial
point of a path, and $\tau(.)$ is the terminal point.

The paths $\mu_i$ and $\nu_i$ are called {\em notches}.  Some or all of
the notches may be trivial.  For $1\leq i< k$, we define $\rho_i$ to be
the path obtained by tightening $\mu_if(\rho_{i-1})\nu_i$.  Since $\rho$
is a closed path, we have $\rho_k=f_\#(\rho_{k-1})$. As before, we call
$\rho_i$ the $i$-th {\em slice} of $\rho$, and the number $k$ is the
duration $\mathcal{D}(\rho)$.

The visible length of $\rho$ is
\[
	\mathcal{V}(\rho)=\mathcal{L}(\rho_0)+\mathcal{L}(\rho_k)+
		\sum_{i=1}^{k-1}\left(\mathcal{L}(\mu_i)+\mathcal{L}(\nu_i)\right).
\]
Finally, we introduce {\em quasi-smooth} hallways: Given some $C\geq 0$,
we say that w hallway $\rho$ is $C$-quasi-smooth if the length of all the
notches is bounded by $C$.

\section{Strata of superlinear growth}\label{nonlinsec}

Throughout this section, let $f\co G\rightarrow G$ be an improved
relative train track map.

In order to track images of edges through the slices of a hallway
$\rho,$ we assign a {\em marking} to each edge.  This assignment will,
in general, involve arbitrary choices, but our arguments will not be
affected by these choices.

\begin{defn}\label{markings}
We begin by marking all edges in the initial slice $\rho_0$ and in all
notches $\mu_i, \nu_i$ with their height.  Assume inductively that all
edges in a slice $\rho_{i-1}$ have been marked, and let $E$ be an edge
of height $r$ in $\rho_{i-1},$ with marking $s.$  Now, consider an edge
$E'$ in $f(E).$  If the height of $E'$ is $r,$ or if $H_s$ is a zero
stratum, then we keep the marking $s.$  If the height of $E'$ is less
than $r,$ then we mark $E'$ by $r.$ This gives us a marking for all
edges in $\mu_if(\rho_{i-1})\nu_i.$

Note that, as we tighten $\mu_if(\rho_{i-1})\nu_i$ to obtain $\rho_i,$
different choices in cancellation (\remref{ambiguity}) may give rise
to different possible markings, but this will not be a problem.

We say that an edge $E$ is marked by a linear/polynomial/exponential
stratum if its marking is $s$ and $H_s$ is
linear/polynomial/exponential.
\end{defn}

The following proposition goes a long way toward proving \thmref{mainthm}.
In fact, if $f$ has no edges of linear growth, then it immediately implies
\thmref{mainthm}.

\begin{prop}\label{nonlinprop}
There exists some constant $K\geq 1$ such that for every hallway $\rho$
and every slice $\rho_i$ of $\rho,$ the number of edges in $\rho_i$ that
are not marked by strata of linear growth is bounded by
$K\mathcal{V}(\rho)$.

Given the improved relative train track map $f\co G\rightarrow G$, the
constant $K$ can be computed.
\end{prop}

In order to streamline the exposition, we will not always make the 
choice of $K$ explicit. However, it will turn out that $K$ can be
chosen to be the product of numbers that can easily be read off
from the train track map.

The intuition of the proof is that once significant growth occurs, it
will be due to the presence of long legal subpaths in exponentially
growing strata or long subsegments of eigenrays of polynomially
growing strata that grow faster than linearly. \lemref{bcc} and
\lemref{sublemma} imply that there is hardly any cancellation
between such subpaths and their surroundings, so that any significant
growth that occurs in a slice will eventually be accounted for by
visible edges.

The following definition will help us understand cancellation in
hallways. For every stratum $H_r$, we define a number $h(H_r)$ in
the following way:
\begin{itemize}
\item If $H_r$ is a constant stratum, then $h(H_r)=0$.
\item If $H_r$ is a nonconstant polynomially growing stratum, i.e.,
$H_r=\{E_r\}$ and $f(E_r)=E_ru_r$, then $h(H_r)$ is the height of $u_r$.
\item If $H_r$ is of exponential growth and $H_{r-1}$ is not a zero stratum,
then $h(E_r)$ is the height of $f(H_r)\cap G_{r-1}$, unless this intersection
does not contain any edges, in which case we let $h(H_r)=\infty$.
\item If $H_r$ is of exponential growth and $H_{r-1}$ is a zero stratum,
then $h(E_r)$ is the height of $f(H_r \cup H_{r-1})\cap G_{r-2}$. We also
let $h(H_{r-1})=h(H_r)$.
\end{itemize}

Essentially, $h(H_r)$ is the index of the highest stratum crossed by
the image of $H_r$, other than $H_r$ itself. We may permute the strata
of $G$ (while preserving the improved train track properties) such that
$h(H_r)>h(H_s)$ implies $r>s$.

Given a stratum $H_s$, we say that the set $S(H_s)=\{H_r|h(H_r)=s\}$ is the
{\em league} of $H_s$, the motivation being that they, in a sense,
``play at the same level.'' If $h(H_r)=\infty$, then $H_r$ does not belong
to any league.

\begin{proof}[Proof of \propref{nonlinprop}]
First of all, we note that if a slice $\rho_i$ has a subpath in a
zero stratum $H_r$, then this subpath is of uniformly bounded length,
and it is surrounded by edges in higher strata (\thmref{ttimproved},
\partref{zerostrat}), so that we have a linear estimate of the number of
edges in $H_r$ in $\rho_i$ in terms of the number of edges in higher
strata.

Let $q$ be the largest (finite) number for which the league $S(q)$ is
nonempty.  Fix some stratum $H_r$ for $r>q$. We want to find a linear
bound on the number of edges in $\rho_i\cap H_r$ in terms of visible
edges. By definition of $S(q)$ and choice of $r$, edges in $\rho_i\cap
H_r$ never cancel with edges from other strata or their images.

If $H_r=\{E_r\}$ is of polynomial growth, then any occurrence of $E_r$
in $\rho_i$ is the image of a visible copy of $E_r$, and $\rho_i$ contains
at most one copy of $E_r$ for each visible copy of $E_r$. Hence, the
number of edges in $\rho_i\cap H_r$ is bounded by the number of visible
edges.

Now, assume that $H_r$ is an exponentially growing stratum.
A slice $\rho_i$ decomposes into $r$-legal subpaths with $r$-illegal turns
in between. By \lemref{bcc}, a subpath whose $r$-length is greater
than $\mathcal{C}_r$ (\eqref{critlen}) will eventually be accounted for
by visible edges since it will not be shortened by cancellation within
slices.

Edges in $H_r$ whose $r$-distance from an illegal turn is less than
$\frac{\mathcal{C}_r}{2}$ may cancel eventually, and $\rho_i$ contains at most
$\mathcal{C}_r$ of them per $r$-illegal turn, so that we only need to find a
bound of the number of $r$-illegal turns in terms of the number of
visible edges. Since the improved train track map $f$ does not create
any $r$-illegal turns, any $r$-illegal turn in $\rho_i$ can be traced back
to a visible illegal turn in $\rho$ (or an illegal turn created by
appending a notch to the image of a slice). This implies that the number
of $r$-illegal turns in $\rho_i$ is bounded by the number of visible
edges in $\rho$.

Summing up, we have bounded the number of edges in
$\rho_i\cap (H_{q+1}\cup H_{q+2} \cup \ldots)$ by a multiple of
the number of visible edges. This establishes the base case of the
proof.

We now assume inductively that the number of edges in $S(p)\cup
S(p+1)\cup \ldots$ has been bounded as a constant multiple of
$\mathcal{V}(\rho).$ We need to find a bound on the number of edges in
$\rho_i\cap H_p$.

We first assume that $H_p=\{E_p\}$ is of polynomial growth. By definition of
$S(p)$, an edge in $\rho_i\cap H_p$ has one of four possible markings:
\begin{itemize}
\item Its marking may be $p,$ indicating that it is the image of a
visible edge, or
\item it may be marked by an exponentially growing stratum in $S(q),$
for some $q\geq p,$ or
\item it may be marked by a superlinear polynomially growing stratum
in $S(q), q\geq p$, or
\item it may be marked by a stratum of linear growth.
\end{itemize}
We are not concerned with edges of the fourth kind.

As before, the number of edges of the first kind in $\rho_i\cap H_p$ is
bounded by the number of visible edges. Let $C$ be the largest number of
copies of $E_p$ that occur in the image of a single edge in an
exponentially growing stratum $H_s$, for $s>p$. Then the number of edges
of the second kind in $\rho_i\cap H_p$ is bounded by $C$ times the
number of exponentially growing edges in $\rho_{i-1}\cap (S(p)\cap
S(p+1)\cap \ldots)$, which in turn is bounded by a multiple of the
number of visible edges.

We have no immediate bound on the number of edges of the third kind.
As we trace the image of such an edge through subsequent slices, one
of three possible events will occur:
\begin{itemize}
\item Either, it eventually maps to a visible edge, or
\item it cancels with an edge of the first or second kind, or
\item it cancels with an edge in the image of a polynomially growing
(possibly linearly growing) edge in $S(p)$.
\end{itemize}
Note that these events may depend on choices in tightening
(\remref{ambiguity}), but once again our estimates will not be
affected by these choices.

The number of edges for which one of the first two events occurs is
clearly bounded by a multiple of the number of visible edges. We only
need to find a bound on the number of edges in an eigenray that eventually
cancel with edges in another eigenray.

\lemref{sublemma} implies that there is a uniform bound on the number of
edges in $H_p$ that cancel when two rays meet, so that we only need to
find a bound on the number of meetings between two rays. Clearly, any
two rays meet at most once.

If an eigenray cancels with segments from more than one
other ray (this is conceivable since a slice may be of the form
$\rho_i=E_rS_1S_2$, where $E_r$ is a polynomially growing edge
in $S(p)$ and $S_1, S_2$ are short segments from rays of edges in
$S(p)$ such that the ray of $E_r$ successively cancels with $S_1$ and
$S_2$), then all except possibly one of these segments cancel completely,
so that they are no longer available for subsequent cancellation. This implies
that the number of meetings of rays is bounded by two times the number of
pieces of rays available for cancellation, which in turn is bounded by the
number of visible edges.

This completes our estimate of the number of edges in $\rho_i\cap H_p$
when $H_p$ is of polynomial growth. We now assume that $H_p$ is of
exponential growth.

The number of subpaths of height $p$ of $\rho_i$ is bounded by the
number of edges of height greater than $p$ in $\rho_i$ plus one.
The contribution of $p$-legal subpaths of $p$-length less than or
equal to $\mathcal{C}_p$ is bounded by $\mathcal{C}_p$ times the number of subpaths of
height $p$, so that we do not need to consider them here. Any
$p$-legal subpaths of length greater than $\mathcal{C}_p$ will eventually
show up in the visible part of $\rho$, so that we do not need to
consider them, either. The remaining edges in $\rho_i\cap H_p$
are at $p$-distance less than $\frac{\mathcal{C}_p}{2}$ from a $p$-illegal
turn. Hence, we only need to find a bound on the number of $p$-illegal
turns in $\rho_i$.

As before, we trace illegal turns in $\rho_i$ back to their origin:
\begin{itemize}
\item An illegal turn may be the image of a visible illegal turn (this
case includes illegal turns created by appending notches $\mu_i, \nu_i$
to the image $f_\#(\rho_{i-1})$ of a slice), or
\item it may come from a illegal turn in the image of an exponentially
growing edge in $S(p)$, or
\item it may be contained in the ray of a polynomially growing edge in
$S(p)$, or
\item it may be contained in a Nielsen path marked by a linear stratum
(\exref{eglinear}).
\end{itemize}
We are not concerned with illegal turns of the fourth type.

The same arguments that we used for polynomially growing $H_p$ yield that
the number of illegal turns of the first and second kind is bounded
by a multiple of the number of visible edges.

Now, let $C$ be the maximum of the number of illegal turns in the images
of polynomially growing edges in $S(p)$. \lemref{fastpoly} yields an
exponent $k_0$ such that for polynomially growing edge $E_r$ in $S(p)$,
$f^{k_0}_\#(u_r)$ contains a long legal segment. This means, in
particular, that if $\rho$ contains a block $f^k_\#(u_r), k\geq k_0$,
then this block contains no more than $C$ illegal turns per long legal
segment. Since long legal segments eventually show up as visible edges,
the number of illegal turns in such blocks is bounded by
$C\mathcal{V}(\rho)$.

The remaining illegal turns are contained in initial subpaths of rays
that contain no more than the first $C+1$ blocks, i.e., there are at
most $C(C+1)$ illegal turns of this kind per ray. Since we already
know that the number of rays is bounded in terms of the number of
visible edges, we are done in this case.

We have now obtained the desired estimate for edges in $\rho_i$
of height $p$ and higher. In particular, this includes all strata
in $S(p-1)$, which completes the inductive step.
\end{proof}

\section{Polynomially growing automorphisms}\label{polysec}

In this section, we establish \thmref{mainthm} in the case of polynomially
growing automorphisms. Specifically, we find estimates for the contribution
of linearly growing edges that we ignored in \propref{nonlinprop}.
As usual, let $f\co G\rightarrow G$ be an improved relative train track map.
Since $f$ is of polynomial growth, every stratum $H_r$ contains only one edge
$E_r$, and we have $f(E_r)=E_r\cdot u_r$, where $u_r$ is some closed path
in $G_{r-1}$. Note that all vertices of $G$ are fixed.

We first record an obvious lemma.

\begin{lemnp}\label{basicnp}
Let $\mu_1, \mu_2$ be Nielsen paths in $G$, and let $\nu$ be some
path in $G$.
\begin{itemize}
\item If $\mu_1$ and $\mu_2$ can be concatenated, then the path
obtained from $\mu_1\mu_2$ by tightening relative endpoints is
also a Nielsen path.
\item If $\mu_1$ and $\nu$ can be concatenated, let $\gamma$ be
the path obtained by tightening $\mu_1\nu$, and let
$\Delta=L(\gamma)-L(\nu)$. Then, for all $k\geq 0$, we have
\[
	L\left(f^k_\#(\gamma)\right)=L\left(f^k_\#(\nu)\right)+\Delta,
\]
and
\[
	-L(\mu_1)\leq \Delta \leq L(\mu_1).
\]
\end{itemize}
\end{lemnp}

We now establish \thmref{mainthm} for automorphisms of linear growth.
This lemma will provide the base case of our inductive proof of
\thmref{mainthm}.

\begin{lem}\label{basecase}
Assume that $f\co G\rightarrow G$ is of linear growth.
If $\rho$ is a smooth hallway, and if $\rho_0$ starts and ends
at vertices, then the lengths of slices of $\rho$ are bounded by
$V(\rho)$, i.e., \thmref{mainthm} holds with $K=1$.
\end{lem}

\begin{proof}
The proof proceeds by induction up through the strata of $G$. The
bottom stratum $H_1$ is constant, so that the lemma trivially holds for
the restriction of $f$ to $H_1$. We now assume that $H_r$ is a linearly
growing stratum, and that the lemma holds for the restriction of $f$
to $G_{p-1}$.

Consider the initial slice $\rho_0$. \remref{howtosplit} yields a
splitting of $\rho_0$ into basic paths of height $p$ and paths in
$G_{p-1}$. The splitting of $\rho_0$ induces a decomposition of
$\rho$ into smooth hallways, so that it suffices to prove the claim
for hallways whose initial slice is a basic path of height $p$ or a
path in $G_{p-1}$.

By induction, we only need to prove the claim if $\rho_0$ is a basic
path of height $p$. If the basic path $\rho_0$ is, in fact, an
exceptional path, then the reasoning of \exref{smoothex} proves
our claim, so that we may assume that $\rho_0$ is not an exceptional
path.

Assume that $\rho_0$ is a basic path of the form $E_p\gamma$. Then, by
\thmref{ttimproved}, \partref{polyim}, there exists some smallest
exponent $m\geq 0$ for which $f^{m+1}_\#(E_p\gamma)$ splits as
$E_p\cdot\gamma'$. Using \remref{howtosplit} once more, we conclude
that $E_p\gamma$ can be expressed as $E_pu_p^{-m}\nu$.

If $D(\rho)\leq m$, then $\rho_0$ $k$-splits as $E_pu_p^{-m}\cdot\nu$.
We can consider the subpaths $E_pu_p^{-m}$ and $\nu$ separately, so that
we are done in this case.

Now assume that $D(\rho)>m$. For $0\leq i\leq m$, we have
$\rho_i=E_pu_p^{i-m}f^i_\#(\nu)$ and
$L(\rho_i)=1+(m-i)L(u_p)+L\left(f^i_\#(\nu)\right)$.
For $m+1\leq i\leq D(\rho)$, we have
$\rho_i=E_pu_p^{i-(m+1)}f^i_\#(u_p\nu)$ and
$L(\rho_i)=1+(i-m-1)L(u_p)+L\left(f^i_\#(\nu)\right)+\Delta$,
where $\Delta$ is defined as in \lemref{basicnp}.

We have
$V(\rho)=2+\left(D(\rho)-1\right)L(u_p)+
L(\nu)+L\left(f^{D(\rho)}_\#(\nu)\right)+\Delta.$
By induction, we have
$L\left(f^i_\#(\nu)\right)\leq L(\nu)+L\left(f^{D(\rho)}_\#(\nu)\right)$
for all $0\leq i \leq D(\rho)$. This immediately implies that
$L(\rho_i)\leq V(\rho)$ for all $0\leq i \leq D(\rho)$.

If $\rho_0=E_p\gamma E_p^{-1}$, we essentially repeat the same
argument. Once more, we can write $\rho_0=E_pu_p^{-m}\nu$, and in
order to use the previous argument, we only need to know that
the lemma holds for $\nu$. This, however, follows from the previous
step, so that we are done.
\end{proof}

We now find estimates on the number of edges emitted by linearly
growing edges, the quantity we ignored in \propref{nonlinprop}.
The idea is to take a hallway and decompose it into smaller and
smaller pieces until all remaining pieces only involve linearly
growing edges and their rays. Simple counting arguments will
give us bounds on the number of the remaining pieces as well
as the lengths of their slices.

\begin{figure} 
\renewcommand{\epsfsize}[2]{0.8\textwidth}
\begin{center}
\psfrag{a}{$\alpha$}
\psfrag{b}{$\beta$}
\psfrag{e}{$E_r$}
\psfrag{u}{$u_r$}
\psfrag{tk}{$t^k$}
\psfrag{cut}{cut}
\psfrag{sawtooth}{sawtooth construction}
\epsfbox{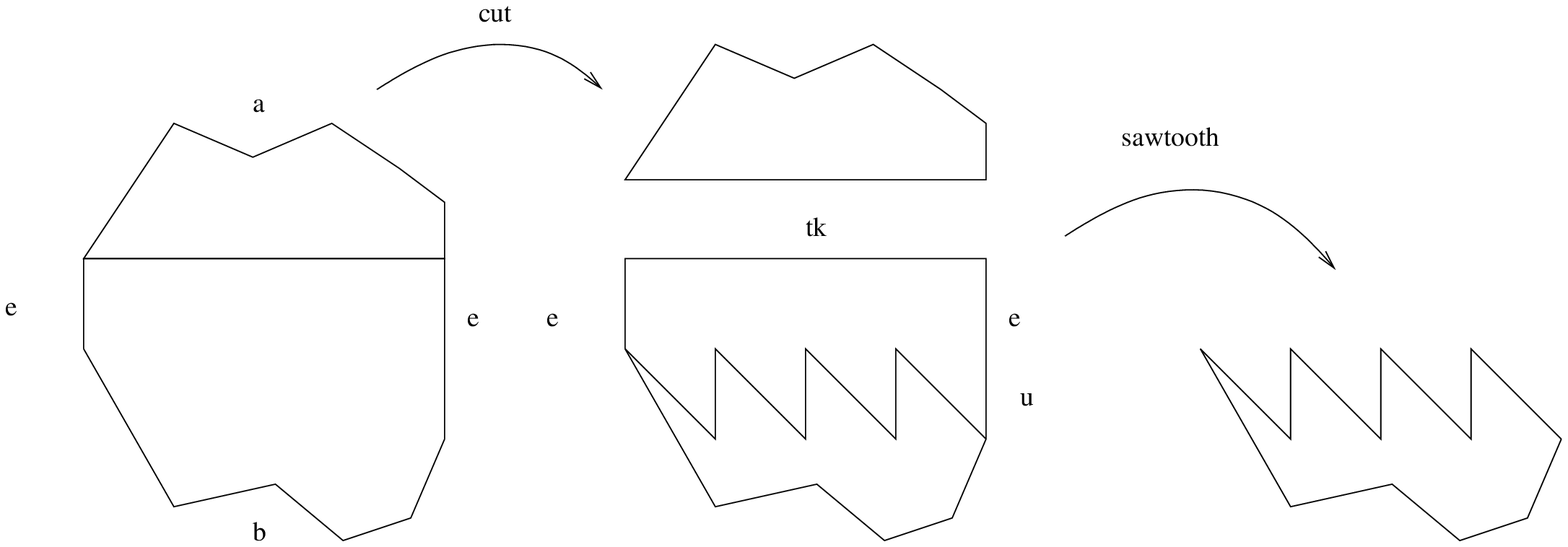}
\end{center}
\caption{Cutting and the sawtooth construction.}
\label{sawtooth}
\end{figure} 

Let $\rho$ be a hallway, and assume that there is a visible edge $E_r$ that
does not cancel within $\rho$, i.e., we can trace its image through the 
slices of $\rho$ until it reappears as another visible edge. Then $\rho$
can be expressed as $\rho=\alpha E_r\beta E_r^{-1}$, and we define two new
hallways $\rho',\rho''$ by tightening $t^{-k}E_r\beta E_r^{-1}$ and
$\alpha t^k$. We say that $\rho'$ and $\rho''$ are obtained from $\rho$ by
{\em cutting} along the trajectory of $E_r$ (\figref{sawtooth}). The
exponent $k$ is the {\em length} of the cut. We say that a hallway $\rho$
is {\em indecomposable} if it does not admit any cuts of length $D(\rho)$.

Now we obtain a new hallway $\sigma$ from $\rho'$ by repeatedly replacing
subwords of the form $t^{-1}E_r$ by $f(E_r)t^{-1}$ and tightening
(\figref{sawtooth}). We refer to this operation as the {\em sawtooth
construction} along the trajectory of $E_r$.

If ${\cal M}$ is a collection of hallways, we let
\[
	V({\cal M})=\sum_{\sigma\in {\cal M}} V(\sigma).
\]

The following lemma lists some basic properties of our two operations.
We say that an edge is of {\em degree} $d$ if $f^k_\#(E)$ grows
polynomially of degree $d$.

\begin{lem}\label{sawtoothprop}
Fix some $C\geq\max\{L(u_r)\}$. Let $\rho$ be a $C$-quasi-smooth hallway in
$G$.  Choose $d>1$ such that the fastest growing edge crossed by $\rho$
grows polynomially of degree $d$.

Obtain a collection $\cal M$ of hallways by cutting along all trajectories
of edges $E$ in $\rho$ of degree $d$.
Let ${\cal M}_1$ be the collection of smooth elements of $\cal M$,
and let ${\cal M}_2$ consist of hallways obtained by performing the
sawtooth construction along all trajectories of $E$ of degree $d$
in those elements of $\cal M$ that are not smooth. Then
\begin{enumerate}
\item The duration of all elements of ${\cal M}_1$ and ${\cal M}_2$ is at
most $D(\rho)$.
\item None of the elements of ${\cal M}_2$ crosses edges of degree $d$,
i.e., they only cross edges of degree at most $d-1$.
\item All elements of ${\cal M}_2$ are $2C$-quasi-smooth.
\item The number of elements of ${\cal M}_2$ is bounded by $2CD(\rho)$.
\item We have
\[
	V({\cal M}_1)+V({\cal M}_2) \leq
		V(\rho)+\left(2CD(\rho)\right)^2.
\]
\end{enumerate}
\end{lem}

\begin{proof}
The first four properties follow immediately from definitions. In order
to prove the fifth property, we just remark that each element of ${\cal M}_2$
has at most $2CD(\rho)$ visible edges that do not appear in $\rho$ itself.
Since ${\cal M}_2$ contains at most $2CD(\rho)$ hallways, the estimate
follows.
\end{proof}

\begin{lem}\label{lowerbound1}
There exists a (computable) constant $C$ with the following property:

Let $\gamma$ be a path of height $r$, starting and ending at vertices,
and assume that $E_r$ is of degree $d>1$. Then, for all $k\geq 0$,
\[
	L(\gamma)+L(f^k_\#(\gamma))\geq Ck^d.
\]
\end{lem}

\begin{proof}
It suffices to prove the lemma if either $\gamma=E_r\gamma'$, or
$\gamma=E_r\gamma'E^{-1}_s$, where $\gamma'$ only involves edges of
degree less than $d$, and $E_s$ is of degree $d$.

In the first case, the claim is obvious. In the second case, we remark
that \lemref{sublemma} guarantees that there is hardly any cancellation
between the rays of $E_r$ and $E_s$, so that the lemma follows.
\end{proof}

The following proposition implies the second part of \thmref{mainthm}
in the case of polynomially growing automorphisms. In particular, it
provides bounds on the number of edges emitted by linearly growing
edges. This is the quantity that we ignored in \propref{nonlinprop}.

\begin{prop}\label{polyprop}
Assume that $f$ represents an automorphism that grows polynomially of degree
$q$. Fix some $C\geq \max\{L(u_r)\}$. There exist computable constants
$K_1\leq K_2 \leq \ldots\leq K_q$ and $K_1'(C),\ldots,K_q'(C)$ such that
\begin{enumerate}
\item If $\rho$ is a smooth hallway whose fastest growing edge is of
degree $d$, and if $\rho_0$ starts and ends at vertices, then
\[
	L(\rho_i)\leq K_dV(\rho)
\]
for all slices $\rho_i$ of $\rho$.
\item If $\rho$ is a $C$-quasi-smooth hallway whose fastest growing edge
is of degree $d$, then in every slice $\rho_i$, the number of edges
emitted by linearly growing edges is bounded by
\[
	K_dV(\rho)+K_d'(C)D^{d+1}(\rho),
\]
so that we have
\[
	L(\rho_i)\leq (K+K_d)V(\rho)+K_d'(C)D^{d+1}(\rho),
\]
where $K$ is the constant from \propref{nonlinprop}.
\end{enumerate}
\end{prop}

\begin{proof}
We prove the proposition by induction on $d$. For $d=1$, the first part
holds with $K_1=1$ because of \lemref{basecase}. Now, assume that $\rho$
is a $C$-quasi-smooth hallway whose fastest growing edge grows of degree
$d=1$. Obtain a collection $\cal M$ of hallways by cutting $\rho$ along
the trajectories of all linearly growing edges that do not cancel within
$\rho$. If $\sigma$ is a smooth element of $\cal M$, then the first
part implies that the number of edges in each $\sigma_i$ emitted by
linearly growing edges is bounded by $V(\sigma)$.

If $\sigma$ is not smooth, then in every slice $\sigma_i$, the number
of edges emitted by linearly growing edges is bounded by
$V(\sigma)+2CD(\rho)$ (It is helpful to keep \exref{bulgeex} in mind).
\lemref{sawtoothprop} yields that $\cal M$ contains no more than
$2CD(\rho)$ pieces that are not smooth. Summing up, we conclude that
every slice of $\rho$ contains at most $V(\rho)+\left(2CD(\rho)\right)^2$
edges emitted by linearly growing edges, so that the second statement
follows with $K_1=1$ and $K_1'(C)=4C^2$.

Now, let $K$ be the constant from \propref{nonlinprop}, and assume
inductively that the proposition holds for some $d\geq 1$. We want
to find some $K_{d+1}$ such that for all hallways $\rho$ whose
fastest growing edge is of degree $d+1$, we have
\[
	L(\rho_i)\leq K_{d+1}V(\rho).
\]
for all slices $\rho_i$. It suffices to prove this with the assumption
that $\rho$ is indecomposable. Then we can perform the sawtooth construction
along all trajectories of edges of degree $d+1$. Since $\rho$ is
indecomposable, we obtain one $C$-quasi-smooth piece $\sigma$ that only
crosses edges of degree $d$ or lower, so that by induction, we conclude
that the number of edges in $\sigma_i$ that were emitted by linearly
growing edges is bounded by
\[
	K_dV(\sigma)+K_d'(C)D^{d+1}(\rho).
\]
We conclude that
\begin{eqnarray*}
	L(\rho_i) &\leq & KV(\rho)+K_dV(\sigma)+K_d'(C)D^{d+1}(\rho) \\
		 &\leq & (K+K_d)V(\rho)+\left(2C+K_d'(C)\right)D^{d+1}(\rho).
\end{eqnarray*}
Using \lemref{lowerbound1}, we can find some constant $M$ such that
\[
	MV(\rho) \geq \left(2C+K_d'(C)\right)D^{d+1}(\rho)
\]
for all indecomposable hallways $\rho$ involving edges of degree $d+1$.
We conclude that the first statement of the proposition holds with
$K_{d+1}=K+K_d+M$.

We now prove the second assertion. Let $\rho$ be a $C$-quasi-smooth hallway.
We obtain two collections ${\cal M}_1, {\cal M}_2$ of hallways by performing
cutting and sawtooth operations as in \lemref{sawtoothprop}.

The elements of ${\cal M}_1$ are smooth hallways, so that for any
$\sigma\in {\cal M}_1$, the previous step yields
\[
	L(\sigma_i)\leq K_{d+1}V(\sigma).
\]

If $\sigma$ is an element of ${\cal M}_2$, then it is a $2C$-quasi-smooth
hallway, and induction yields that in every slice of $\sigma$, the number
of edges emitted by linearly growing edges is bounded by
\[
	K_dV(\sigma)+K_d'(2C)D^{d+1}(\sigma).
\]

Summing over all elements of ${\cal M}_1$ and ${\cal M}_2$, we conclude that
every slice $\rho_i$ of $\rho$ contains at most
\begin{eqnarray*}
K_{d+1}V({\cal M}_1)+K_dV({\cal M}_2)+2CD(\rho)\cdot K_d'(2C)D^{d+1}(\rho)\\
\leq K_{d+1}V(\rho)+4C^2K_dD^2(\rho)+2CK_d'(2C)D^{d+2}(\rho)
\end{eqnarray*}
edges emitted by linearly growing edges, so that the second statement of
the proposition holds with
\[
	K_{d+1}'(C)=4C^2K_d+2CK_d'(2C).
\]
\end{proof}

\begin{rem}
The estimates of \propref{polyprop} are rather crude; lots of edges are
counted several times rather than just once. I opted to present
the most straightforward estimates rather than tightest ones.
\end{rem}

\section{Proof of the main result}\label{pfsec}

We now extend the techniques and results of \propref{polysec} to
arbitrary automorphisms. The presence of exponentially growing strata
will turn out to be a mixed blessing. On the one hand, they make for
rather simple counting arguments as polynomial contributions as in
\propref{polyprop} are easily dwarfed by exponential growth. On the
other hand, we will need to consider more complicated decompositions
of hallways.

As usual, let $f\co G\rightarrow G$ be an improved relative train track
map.  Any statements regarding the computability of constants assume
that we are given such a map.  After permuting the strata as necessary,
we may assume that if $H_r$ and $H_s$ are truly polynomial strata and
$r>s$, then the degree of $H_r$ is at least as large as that of $H_s$.
Throughout this section, let $K$ be the constant from
\propref{nonlinprop}.

If $H_r$ is an exponentially growing stratum, then we fix some
$L>\mathcal{C}_r$, and we replace $f$ by $f^M$, where $M$ is the exponent
from \lemref{tricho} for this choice of $L$. After replacing $f$
by a power yet again if necessary, we may assume that the image
of each edge in $H_r$ contains at least $L$ edges in $H_r$.
If $H_r$ supports a closed Nielsen path $\tau$, then the initial
and terminal edges of $\tau$ are partial edges in $H_r$, and we
may assume that the image of each of them also contains at least
$L$ edges in $H_r$. We say that a legal path of height $r$
is {\em long} if it contains at least $L$ edges in $H_r$.

We first record an exponential version of \lemref{lowerbound1}.

\begin{lem}\label{lowerbound2}
Let $H_r$ be an exponentially growing stratum or a fast polynomial
stratum.  Then there exists a computable constant $\lambda>1$ such that
if $\sigma$ is a circuit in $G_r$ or a path starting and ending at fixed
vertices, then either $\sigma$ is a concatenation of Nielsen paths of
height $r$ and subpaths in $G_{r-1}$, or we have
\[
	\mathcal{L}(\sigma)+\mathcal{L}(f^k_\#(\sigma))\geq \lambda^k
\]
for all $k\geq 0$.
\end{lem}

\begin{proof}
If $H_r$ is an exponentially growing stratum, we need to distinguish two
cases: First, assume that for some $i\geq 0$, $f^i_\#(\sigma)$ is a
concatenation of Nielsen paths and subpaths in $G_{r-1}$. Since $\sigma$
starts and ends a fixed vertices, we conclude that $\sigma$ itself is a
concatenation of Nielsen paths and subpaths in $G_{r-1}$, so that there
is nothing to show in this case.

Let $\lambda_-,N_0$ be the constants from \lemref{expdecay},
and assume that for all $i\geq 0$, $f^i_\#(\sigma)$ is not a
concatenation of Nielsen paths and subpaths in $G_{r-1}$.
Let $i_0$ be the smallest index for which $f^{i_0}_\#(\sigma)$
contains a long legal segment. Then, using \lemref{expdecay} and
\lemref{tricho}, we see that $\mathcal{L}(\sigma)\geq \lambda_-^{i_0}$.
Moreover, we have $\mathcal{L}(f^k_\#(\sigma))\geq \lambda_r^{k-i_0}$.

If we let $\lambda=\sqrt{\min\{\lambda_-,\lambda_r\}}$, then we have
$\lambda_-^{i_0}+\lambda_r^{k-i_0}\geq \lambda^k$. Hence, we have
$\mathcal{L}(\sigma)+\mathcal{L}(f^k_\#(\sigma))\geq
\lambda_-^{i_0}+\lambda_r^{k-i_0} \geq \lambda^k$.

If $H_r=\{E_r\}$ is a fast polynomial stratum, then we argue similarly,
using \lemref{fastpoly} and \thmref{ttimproved}, \partref{polyim}.
\end{proof}

If $H_r$ is an exponentially growing stratum, we let $T_r$ equal the
length of the longest path in $f(H_r)\cap G_{r-1}$.
We fix another constant $S_r>0$ with the following property: Let $\gamma$
be a path in $G_{r-1}$. If $\mathcal{L}(\gamma)\geq S_r$, then
$\mathcal{L}(f_\#(\gamma))>3T_r$ and $\mathcal{L}(f^2_\#(\gamma))>3T_r$, and
if $\mathcal{L}(\gamma)\leq T_r$, then $\mathcal{L}(f_\#(\gamma))<S_r$ and
$\mathcal{L}(f^2_\#(\gamma))<S_r$. We can easily compute a suitable
value $S_r$ given the train track map $f$.
We say that a path $\gamma$ in $G_{r-1}$ is {\em $r$-significant} if 
$\mathcal{L}(\gamma)\geq S_r$.

If $H_r$ is an exponentially growing stratum, and $\rho$ is a
$C$-quasi-smooth hallway of height $r$, then we need to develop an
understanding of the lengths of components of $\rho_i\cap G_{r-1}$,
i.e., we need to study subpaths in $G_{r-1}$.  Intuitively, we will
accomplish this by carving out subhallways in $G_{r-1}$.

Consider a maximal subpath $\gamma\subset G_{r-1}$ of some slice
$\rho_a,$ i.e., $\rho_a$ can be expressed as $\alpha\gamma\beta,$ and
$\alpha$ (resp.\ $\beta$) is either trivial or ends (resp.\ starts) with
a (possibly partial) edge in $H_r.$  We begin the construction of a new
hallway $\rho'$ by letting $\rho_0'=\gamma$.

Now, assume inductively that we have defined the slice $\rho_{i-1}'$
such that $\rho_{i-1}'$ is a maximal subpath of $\rho_{a+i-1}$ in
$G_{r-1}$ (we write $\rho_{a+i-1}=\alpha\rho_{i-1}'\beta$), and
recall that the slice $\rho_{a+i}$ is obtained by tightening
$\mu_{a+i}f(\alpha\rho_{i-1}'\beta)\nu_{a+i}$.  We define the notch $\mu_i'$
by taking the maximal terminal subpath in $G_{r-1}$ of the path obtained
from $\mu_{a+i}f(\alpha)$ by tightening.  Similarly, we define the notch
$\nu_1'$ by tightening the maximal initial subpath in $G_{r-1}$ of the
path obtained from $f(\beta)\nu_{a+i}$ by tightening.  Observe that
tightening $\mu_i'\rho_{i-1}'\nu_i'$ yields a maximal subpath in $G_{r-1}$
of $\rho_{a+i}$, and that the length of $\mu_i'$ and $\nu_i'$ is bounded
by $C+T_r$.  We iterate this procedure until we reach a point where
tightening $\mu_{i+1}'f(\rho_i')\nu_{i+1}'$ yields a trivial path.

By applying this construction wherever possible, we obtain a fan of
$C+T_r$-quasi-smooth hallways in $G_{r-1}$.  Let $\mathcal{M}$ be the
set of maximal elements of this fan.  We let ${\cal M}_1$ be the
collection of smooth hallways in $\mathcal{M}$, and we let ${\cal M}'_2$
be the collection of hallways in $\mathcal{M}$ that are not smooth.

Let $\sigma$ be an element of ${\cal M}'_2$, and assume that there
exists some $0<i<\mathcal{D}(\sigma)$ such that
$\mathcal{L}(\sigma_i)<S_r$.  Then we obtain two new hallways
$\sigma',\sigma''$ from $\sigma$ by letting $\sigma'_j=\sigma_j$ for
$0\leq j\leq i$ and $\sigma''_j=\sigma_{i+j}$ for $0\leq j\leq
\mathcal{D}(\sigma)-i$; we may think of this operation as cutting
$\sigma$ along $\sigma_i$. We obtain a collection of hallways ${\cal
M}_2$ by performing all possible cuts of this kind on all elements of
${\cal M}'_2$.

If $\sigma\in {\cal M}_1\cup{\cal M}_2$, we say that $\sigma$ {\em intersects}
a slice $\rho_i$ if one of the slices of $\sigma$ is a subpath of $\rho_i$.
When looking for bounds on the lengths of a slice $\rho_i$, we need to
find bounds on the lengths of slices of hallways $\sigma$ that intersect
$\rho_i$.

\def\condA#1{Condition~$A_{#1}$}
Fix some stratum $H_r$. We say that the map $f$ satisfies {\em
\condA{r}} if for any $C\geq 0$, there exist computable constants $K_r$,
$K_r'(C)$, and an exponent $d\geq 1$, such that the following two
conditions hold:
\begin{itemize}
\item If $\rho$ is a smooth hallway in $G_r$ such that the slice
$\rho_0$ starts and ends at fixed vertices, then
\[
	\mathcal{L}(\rho_i)\leq K_r\mathcal{V}(\rho)
\]
for all slices $\rho_i$.
\item If $\rho$ is a $C$-quasi-smooth hallway in $G_r$, then
\[
	\mathcal{L}(\rho_i)\leq K_r\mathcal{V}(\rho)+K'_r(C)\mathcal{D}^d(\rho).
\]
\end{itemize}

If $H_r$ is an exponentially growing stratum, then a hallway of height $r$
is {\em admissible} if all its slices start and end at fixed vertices or
at points in $H_r$.

\begin{lem}\label{subhallwaylem}
Let $H_r$ be an exponentially growing stratum, and assume that \condA{r-1}
holds. Then, given some $C\geq 0$, there exist computable constants
$C_1,C_2\geq 1$ with the following property: If $\rho$ is an admissible
$C$-quasi-smooth hallway of height $r$,then
\[
	\mathcal{L}(\rho_i)\leq C_1\mathcal{V}(\rho)+
		C_2\sum_{{\sigma\in{\cal M}_2\atop\textrm{$\sigma$ intersects $\rho_i$}}
		\atop\textrm{in an $r$-significant segment}}\mathcal{D}^d(\sigma)
\]
for every slice $\rho_i$ of $\rho$.
\end{lem}

\begin{proof}
Since $\rho$ is admissible, all slices of $\sigma\in {\cal M}_1$ start
and end at fixed vertices unless $\sigma_0$ is contain in a zero stratum,
in which case all slices $\sigma_i$ for $i>0$ start and end at fixed
vertices. Moreover, if $\sigma_0$ is contained in a zero stratum, then
$\mathcal{L}(\sigma_1)=\mathcal{L}(\sigma_0)$. By \condA{r-1}, we have
\[
	\mathcal{L}(\sigma_i)\leq K_{r-1}\mathcal{V}(\sigma)
\]
for all slices $\sigma_i$ of $\sigma\in {\cal M}_1$.

Fix some slice $\rho_i$ of $\rho$. Using \propref{nonlinprop} and
\condA{r-1}, we see that
\begin{eqnarray*}
	\mathcal{L}(\rho_i) &\leq & K\mathcal{V}(\rho)+
		\sum_{\sigma\in{\cal M}_1\atop\textrm{$\sigma$ intersects $\rho_i$}}
			K_{r-1}\mathcal{V}(\sigma) \\
	& & +\sum_{\sigma\in{\cal M}_2\atop\textrm{$\sigma$ intersects $\rho_i$}}
			\left(K_{r-1}\mathcal{V}(\sigma)+K'_{r-1}(C+T_r)\mathcal{D}^d(\sigma)\right).\\
\end{eqnarray*}

Consider some $\sigma\in {\cal M}_1$ that intersects $\rho_i$. If the
initial slice of $\sigma$ is not visible in $\rho$, then, as we noted
before, its length is bounded by $T_r$. Similarly, if the terminal slice
of $\sigma$ is not visible in $\rho$, then its length is also bounded by
$T_r$. The number of elements of ${\cal M}_1$ that intersect $\rho_i$
is bounded by $K\mathcal{V}(\rho)$. Putting it all together, we conclude that
\[
  \sum_{\sigma\in{\cal M}_1\atop\textrm{$\sigma$ intersects $\rho_i$}}\mathcal{V}(\sigma)
		\leq (2KT_r+1)\mathcal{V}(\rho).
\]
Similarly, using the fact that elements of ${\cal M}_2$ are
$C+T_r$-quasi-smooth, and that their initial and terminal slices are either
visible in $\rho$ or of length less than $S_r$, we see that
\[
  \sum_{\sigma\in{\cal M}_2\atop\textrm{$\sigma$ intersects $\rho_i$}}\mathcal{V}(\sigma)
	\leq (2KS_r+1)\mathcal{V}(\rho)+
	2(C+T_r)\sum_{\sigma\in{\cal M}_2\atop
		\textrm{$\sigma$ intersects $\rho_i$}}\mathcal{D}(\sigma).
\]
Since $\rho_i$ contains at most $K\mathcal{V}(\rho)$ subpaths in $G_{r-1}$, the
total contribution of subpaths in $G_{r-1}$ that are not $r$-significant is
bounded by $KS_r\mathcal{V}(\rho)$. Letting $C_1=K+2K_{r-1}(K(S_r+T_r)+1)+KS_r$ and
$C_2=K'_{r-1}(C+T_r)+2(C+T_r)$, we conclude that
\[
	\mathcal{L}(\rho_i) \leq C_1\mathcal{V}(\rho)+
		C_2\sum_{{\sigma\in{\cal M}_2\atop\textrm{$\sigma$ intersects $\rho_i$}}
			\atop\textrm{in an $r$-significant segment}}
			\mathcal{D}^d(\sigma).
\]
\end{proof}

\lemref{subhallwaylem} shows that from now on, we may focus on the
polynomial contribution of nonsmooth hallways in $G_{r-1}$ that
intersect a given slice $\rho_i$ in an $r$-significant subpath. In
particular, if the initial slice $\rho_0$ happens to be an $r$-legal
path, then
\[
	\mathcal{L}(\rho_i)\leq C_1\mathcal{V}(\rho)
\]
for all slices $\rho_i$ since ${\cal M}_2$ is empty in this case.

\begin{lem}\label{Cqs1}
Let $H_r$ be an exponentially growing stratum, and assume that
\condA{r-1} holds. Given some $C>0$, there exist computable
constants $C_1,C_2$ with the following property:
If $\rho$ is an admissible $C$-quasi-smooth hallway of height $r$, such
that for every slice $\rho_i$ except possibly the last one, $f_\#(\rho_i)$
does not contain a legal segment of length at least $L$, then
\[
	\mathcal{L}(\rho_i)\leq C_1\mathcal{V}(\rho)+C_2\mathcal{D}^{d+1}(\rho)
\]
for all slices $\rho_i$.
\end{lem}

\begin{proof}
By \lemref{subhallwaylem}, we may restrict our attention to elements of
${\cal M}_2$ that intersect a given slice $\rho_i$ in an $r$-significant
subpath.  Let
\[
  D=\sum_{{\sigma\in{\cal M}_2\atop\textrm{$\sigma$ intersects $\rho_i$}}
	\atop\textrm{in an $r$-significant segment}} \mathcal{D}^d(\sigma).
\]

We first claim that the number of $r$-significant subpaths in $G_{r-1}$ in a
slice $\rho_i$ is bounded by $N(\rho_i)$. By choice of $S_r$, an
$r$-significant subpath in $G_{r-1}$ will not cancel completely when
$f(\rho_i)$ is tightened to $f_\#(\rho_i)$.

If there were two such subpaths in one legal segment of $\rho_i$, then
there would be a legal segment in $H_r$ in between. Since we assumed that
$\mathcal{L}(f(E)\cap H_r)\geq L$ for each edge in $H_r$, the $r$-length of
the image of this legal segment is at least $L$, which means that
the slice $\rho_{i+1}$ contains a legal segment of length at least
$L$, contradicting our assumption. This proves the claim if $H_r$
does not support a closed Nielsen path, as in this case, the number
of legal segments in $\rho_i$ equals $N(\rho_i)$.

If $H_r$ supports a closed Nielsen path, then a legal segment of
$\rho_i$ that is adjacent to an illegal turn contained in a Nielsen
subpath of $\rho_i$ cannot contain an $r$-significant subpath in $G_{r-1}$.
If such a segment contained an $r$-significant subpath in $G_{r-1}$, then
$f_\#(\rho_i)$ would contain a legal segment of $r$-length $L$ because
both the initial and terminal partial edge of the Nielsen path of $H_r$
map to legal segments of $r$-length at least $L$. This implies that the
number of $r$-significant subpaths in $G_{r-1}$ is bounded by $N(\rho_i)$.

Now, fix some slice $\rho_i$. We make the worst-case assumption that
every legal segment of $\rho$ that is not adjacent to an illegal turn
contained in a Nielsen subpath contains an $r$-significant subpath in $G_{r-1}$
that is a slice of a hallway $\sigma\in {\cal M}_2$ of duration $j\geq i$.
The number of such hallways whose duration is a given number $j\geq i$
is bounded by $N(\rho_j)+1$. We conclude that
\[
	D \leq \sum_{j=i}^{\mathcal{D}(\rho)}N(\rho_j)j^d.
\]

Choosing $\lambda$ according to \lemref{expdecay},
we conclude that $N(\rho_{i+1})\leq \lambda^{-1}N(\rho_i)+1+2C$, as
$\rho$ is $C$-quasi-smooth. This implies, inductively, that
\[
	N(\rho_i)\leq \lambda^{-i}N(\rho_0)+2(1+C)\sum_{j=0}^{i-1}\lambda^{-j}
		\leq \lambda^{-i}N(\rho_0)+\frac{\lambda}{\lambda-1}(1+2C).
\]
We choose some $B\geq \sum_{j=0}^{\infty}\lambda^{-j}j^d$, and we conclude
that
\[
	D\leq \sum_{j=0}^{\mathcal{D}(\rho)}N(\rho_j)j^d\leq
		B\mathcal{V}(\rho)+\frac{\lambda}{\lambda-1}(1+2C)\mathcal{D}^{d+1}(\rho),
\]
since $N(\rho_0)\leq \mathcal{V}(\rho)$.

If $C'_1,C'_2$ are the constants from \lemref{subhallwaylem}, then the
lemma holds with $C_1=C'_1+C'_2B$ and
$C_2=\frac{\lambda}{\lambda-1}(1+2C)C'_2$.
\end{proof}

\def\propB{Property~$B$}
Let $H_r$ be an exponentially growing stratum, and let $N_0$ be the
constant from \lemref{expdecay}. We say that an admissible smooth
hallway $\rho$ of height $r$ has {\em \propB} if for all slices
$\rho_i$, $\rho_i$ contains no long $r$-legal segment, or
$N(\rho_{i-1})<N_0.$

\begin{lem}\label{latebloom}
Let $H_r$ be an exponentially growing stratum, and assume that \condA{r-1}
holds. Let $N_0$ be the constant from \lemref{expdecay}. There exist
computable constants $C_1,C_2$ with the following property: If $\rho$ is
an admissible smooth hallway of height $r$ that satisfies \propB, then
\[
	\mathcal{L}(\rho_i)\leq C_1\mathcal{V}(\rho)+C_2\mathcal{D}^{d+1}(\rho)
\]
for all slices $\rho_i$.
\end{lem}

\begin{proof}
If no slice of $\rho$ contains a long legal segment, then the claim follows
from \lemref{Cqs1}. Otherwise, let $i_0$ be the smallest index for which
$\rho_{i_0}$ contains a long legal segment. By choice of $i_0$, $\rho_{i_0-1}$
does not contain a long legal segment, and by hypothesis, we have
$N(\rho_{i_0-1})<N_0$. If $i<i_0$, then, choosing $D$ as in the proof
of \lemref{Cqs1}, we conclude that
\begin{eqnarray*}
	D & \leq &
    \left(\sum_{j=0}^{i_0-1}N(\rho_j)j^d\right)+N_0\mathcal{D}^d(\rho) \\
	  & \leq & B\mathcal{V}(\rho)+N_0\mathcal{D}^{d+1}(\rho),
\end{eqnarray*}
so that the lemma holds for all $\rho_i$ with $i<i_0$.

For $i\geq i_0$, $\rho_i$ splits as a concatenation of long $r$-legal
paths and subpaths that contain illegal turns and no long legal subpaths.
Each slice may, conceivably, contain slices of $N(\rho_{i_0-1})<N_0$ 
hallways of duration $\mathcal{D}(\rho)$. The polynomial contribution of these
hallways is bounded by $N_0\mathcal{D}^d(\rho)$.

In addition, the number of short legal segments around illegal turns is
at most $2N_0$. Each of them contains not more than one $r$-significant subpath
in $G_{r-1}$, belonging to a hallway of duration at most
$\mathcal{D}(\rho)-i_0$.  The polynomial contribution of these paths is
bounded by $2N_0(\mathcal{D}(\rho)-i_0)^d$.

Now, since $\rho_{i_0}$ contains a long legal segment, the length of
$\rho_{\mathcal{D}(\rho)}=f^{\mathcal{D}(\rho)-i_0}_\#(\rho_{i_0})$ is at
least $\lambda_r^{\mathcal{D}(\rho)-i_0}$. We can easily find some $B'>0$
such that $B'\lambda_r^k\geq 2N_0k^d$ for all $k\geq 0$. We conclude that for
the sum of all polynomial contributions in $\rho_i$, we have
\[
	2N_0(\mathcal{D}(\rho)-i_0)^d+
        N_0\mathcal{D}^d(\rho)\leq B'\mathcal{V}(\rho)+N_0\mathcal{D}^d(\rho),
\]
which completes the proof of the lemma.
\end{proof}

The remaining two lemmas deal with arbitrary smooth hallways of height
$r$ as well as quasi-smooth hallways by essentially decomposing them
into pieces of the kind that we analyzed in the previous lemmas.

\begin{lem}\label{smoothlem}
Let $H_r$ be an exponentially growing stratum, and assume that
\condA{r-1} holds. Then there exist computable constants $C_1,C_2$
with the following property: If $\rho$ is a smooth admissible hallway
of height $r$, then
\[
	\mathcal{L}(\rho_i)\leq C_1\mathcal{V}(\rho)+C_2\mathcal{D}^{d+1}(\rho)
\]
for all slices $\rho_i$.
\end{lem}

\begin{proof}
Let $\lambda_-,N_0$ be the constants from \lemref{expdecay}. As in
the proof of \lemref{lowerbound2}, we let
$\lambda=\sqrt{\min\{\lambda_-,\lambda_r\}}$, and we remark that for
$0\leq j\leq k$, we have $\lambda^j_-+\lambda^{k-j}_r\geq \lambda^k$.
This basic estimate will be crucial in the proof of this lemma. We
choose some $B>0$ such that $B\lambda^k>k^{d+1}$ for all $k\geq 0$.

Let $C'_1,C'_2$ be the maximum of the corresponding constants from the
previous lemmas. We will see that the lemma holds with $C_1=C'_1+3BC'_2$
and $C_2=C'_2$.

We first observe that if $\rho$ satisfies \propB, then the lemma
follows from \lemref{latebloom}.  If $\rho_0$ contains long legal
segments, we can split $\rho_0$ into long $r$-legal subpaths and
neighborhoods of illegal turns (i.e., illegal turns surrounded by legal
paths whose length is at most $\frac{C_r}{2}$).  Split $\rho_0$ as
$\rho_0=\alpha_{0;1}\beta_{0;1}\alpha_{0;2}\cdots\alpha_{0;m}\beta_{0;m}$,
where all subpaths $\alpha_{0;i}$ are long legal segments, and all
subpaths $\beta_{0;i}$ are neighborhoods of illegal turns. Such a
decomposition of $\rho_0$ induces a decomposition of $\rho$ into
hallways, and we can choose the decomposition such that all resulting
pieces are admissible, and that the legal segments are as long as
possible, subject to admissibility.  We write
$\alpha_{j;i}=f^j_\#(\alpha_{0;i})$ and
$\beta_{j;i}=f^j_\#(\beta_{0;i})$.

Let $k=\mathcal{D}(\rho)$. For each long legal subpath $\alpha_{0;i}$,
\lemref{subhallwaylem} yields that
$\mathcal{L}(\alpha_{j;i})\leq
C_1(\mathcal{L}(\alpha_{0;i})+\mathcal{L}(\alpha_{k;i}))$,
for all $0\leq j\leq k$. Since $\alpha_{0;i}$ is a long legal segment,
we have $\mathcal{L}(\alpha_{k;i})\geq \lambda_r^k\geq \lambda^k$.

If the hallway defined by $\beta_{0;i}$ satisfies \propB, then we have
$\mathcal{L}(\beta_{j;i})\leq
C'_1(\mathcal{L}(\beta_{0;i})+\mathcal{L}(\beta_{k;i}))+C'_2k^{d+1}$,
and we have $k^{d+1}\leq B\mathcal{L}(\alpha_{k;i})$, hence
\[
	\mathcal{L}(\beta_{j;i})\leq
    C'_1(\mathcal{L}(\beta_{0;i})+\mathcal{L}(\beta_{k;i}))
		+BC'_2\mathcal{L}(\alpha_{k;i}),
\]
i.e., we can find a legal segment adjacent to $\beta_{0;i}$ whose
contribution to the visible edges of $\rho$ dominates the possible
polynomial contribution of $\beta_{0;i}$. This takes care of the long
legal segments in $\rho_0$ as well as the subpaths that satisfy \propB.
Hence, we only need to deal with those paths that do not satisfy \propB.
Assume that for some $0\leq i\leq m$, $\beta_{0;i}$ is one of
them.

Then there exists some $j_0$ such that $\beta_{j_0;i}$
contains a long legal segment, but $\beta_{j_0-1;i}$ does not,
and $N(\beta_{j_0-1;i})\geq N_0$.

As before, we split $\beta_{j_0;i}$ into long legal segments and
neighborhoods of illegal turns, obtaining a decomposition
$\beta_{j_0,i}=\alpha_{j_0;i,0}\beta_{j_0;i,0}\cdots
\alpha_{j_0;i,m}\beta_{j_0;i,m}$, where $\alpha_{j_0;i,k}$ are $r$-legal
subpaths, and $\beta_{j_0;i,k}$ are neighborhoods of illegal turns. We
can find splittings
$\beta_{j;i}=\alpha_{j;i,0}\beta_{j;i,0}\cdots\alpha_{j;i,m}\beta_{j;i,m}$
for all $0\leq j\leq k$, such that $f_\#(\alpha_{j;i,k})=\alpha_{j+1;i,k}$
and $f_\#(\beta_{j;i,k})=\beta_{j+1;i,k}$. We may choose those splitting
such that the resulting pieces are admissible, and such that the
legal segments $\alpha_{j_0;i,k}$ are as long as possible, subject to
admissibility.

Now, fix on one subpath $\alpha_{j_0;i,k}$. If $N$ is the number of
$r$-significant subpaths in $G_{r-1}$ in $\alpha_{j_0;i,k}$, then
$\alpha_{j_0-1;i,k}$ contains at least $N$ legal segments containing
$r$-significant subpaths in $G_{r-1}$. By \lemref{expdecay}, we have
$\mathcal{L}(\beta_{0;i})\geq N(\beta_{j_0;i})\geq
\lambda_-^{j_0-1}N(\beta_{j_0-1;i})$, so that we can find
$\lambda_-^{j_0-1}N$ illegal turns in $\beta_{0;i}$, and we can find
$\lambda_r^{k-j_0}$ edges in $\beta_{k,i}$. Using our earlier estimate,
we see that $(\lambda_-^{j_0-1}+\lambda_r^{k-j_0})N\geq
\lambda_-^{-1}\lambda^kN$.

The polynomial contribution of the $r$-significant subpaths in $G_{r-1}$ of
$\alpha_{j_0;i,k}$ is bounded by $K'_{r-1}(T_r)Nk^{d+1}\leq
BK'_{r-1}(T_r)N\lambda^k$, i.e., it is dominated by corresponding
visible edges.

This leaves us to deal with the adjacent subpaths $\beta_{j_0;i,k}$
and $\beta_{j_0;i,k-1}$. If $\beta_{0;i,k}$ satisfies \propB,
then its polynomial contribution is bounded by $C'_2k^{d+1}$, which in
turn is bounded by $BC'_2\lambda^k$.

This takes care of the legal segments $\alpha_{j_0;i,k}$ as well as
those neighborhoods of illegal turns that satisfy \propB.  We apply
the previous reasoning to the remaining paths $\beta_{j_0;i,k}$,
completing the proof of the lemma.
\end{proof}

\begin{lem}\label{Cqs2}
Let $H_r$ be an exponentially growing stratum, and assume that
\condA{r-1} holds. Given some $C>0$, there exist computable
constants $C_1,C_2$ with the following property:
If $\rho$ is an admissible $C$-quasi-smooth hallway of height $r$,
then
\[
	\mathcal{L}(\rho_i)\leq C_1\mathcal{V}(\rho)+C_2\mathcal{D}^{d+3}(\rho)
\]
for all slices $\rho_i$.
\end{lem}

\begin{proof}
The idea of this proof is to decompose the hallway $\rho$ into
pieces that are either smooth or $C$-quasi-smooth satisfying the
hypothesis of \lemref{Cqs1}.

In order to find this decomposition, we introduce {\em trajectories}
of points in $H_r$.  This definition may be affected by the choices made
when tightening (\remref{ambiguity}).  In order to avoid ambiguities,
for each index $1\leq i<D(\rho)$, we fix a sequence of elementary
cancellations that turn $\mu_i\rho_{i-1}\nu_i$ into $\rho_i$.

If $p$ is a point in $\rho_i\cap H_r$, we consider
its image $f(p)$ in $f(\rho_i)$. We say that $p$ {\em survives}
if $f(p)$ is contained in $H_r$ and if $f(p)$ is not contained in an edge
that cancels when $f(\rho_i)$ is tightened to $f_\#(\rho_i)$. If $p$
survives, then $f(p)$ is contained in $\rho_{i+1}$, or it is contained
in the parts of $f_\#(\rho_i)$ that cancel when
$\mu_{i+1}f_\#(\rho_i)\nu_{i+1}$ is tightened to $\rho_{i+1}$.

Thinking of the hallway $\rho$ as spanning a (possibly singular) disk,
we draw a line segment (in this disk) from the surviving points in each
slice to their images. If $p$ is a point in a visible edge such that $p$
and all its images survive, then $p$ defines a line starting and ending
in visible edges, called the {\em trajectory} of $p$. The trajectories
of two points need not be disjoint, but that does not concern us here.

We say that two trajectories are {\em parallel} if their initial points
are both contained in $\rho_0$ or both contained in the same notch, and
if their terminal points are both contained in
$\rho_{\mathcal{D}(\rho)}$ or both contained in the same notch. The
crucial observation is that equivalence classes of parallel trajectories
are closed subsets of the disk spanned by $\rho$, so that in every
equivalence class, we can find trajectories of two points $p_1,p_2$ that
are extremal in the following sense: If $p$ is a point whose trajectory
is parallel to those of $p_1$ and $p_2$, then $p$ is located between
$p_1$ and $p_2$.

We now cut $\rho$ along the extremal trajectories of all equivalence
classes of parallel trajectories, obtaining pieces that are either smooth
or $C$-quasi-smooth. Moreover, all the resulting pieces are admissible.
Let ${\cal M}_1$ be the collection of smooth pieces and ${\cal M}_2$ the
collection of pieces that are not smooth. Note that
$\mathcal{V}({\cal M}_1)+\mathcal{V}({\cal M}_2)=\mathcal{V}(\rho)$.

We now claim that all elements of ${\cal M}_2$ satisfy the hypothesis
of \lemref{Cqs2}. Suppose otherwise, i.e., there exists some
$\sigma\in {\cal M}_2$ such that for some slice $\sigma_i$,
$f_\#(\sigma_i)$ contains a legal segment of length at least $L$.
Within the interior of this legal segment, we can find some point $p$
such that all images of $p$ survive in subsequent slices. Since $p$ is
the image of surviving points, we obtain a trajectory along which we
can cut $\sigma$, contradicting the fact that we obtained $\sigma$ by
cutting $\rho$ along extremal trajectories.

By \lemref{smoothlem}, there are constants $C'_1,C'_2$ such that for
every $\sigma\in {\cal M}_1$ and every slice $\sigma_i$ of $\sigma$,
we have
\[
	\mathcal{L}(\sigma_i)\leq C'_1\mathcal{V}(\sigma)+C'_2\mathcal{D}^{d+1}(\sigma),
\]
and by \lemref{Cqs2}, there are constants $C''_1,C''_2$ such that
\[
	\mathcal{L}(\sigma_i)\leq C''_1\mathcal{V}(\sigma)+C''_2\mathcal{D}^{d+1}(\sigma)
\]
for every slice $\sigma_i$ of every $\sigma\in {\cal M}_2$.

There are at most $2(\mathcal{D}(\rho)-1)$ notches, so that the number of
equivalence classes of parallel trajectories is bounded by
$(2(\mathcal{D}(\rho)-1)+1)^2$ (another extremely crude estimate, but it'll
do).  Since we cut along no more than two trajectories per equivalence class,
we obtain no more than
\[
    2\left(2\left(\mathcal{D}(\rho)-1\right)+1\right)^2+1\leq
    8\mathcal{D}^2(\rho)
\]
pieces.  Letting $C_1=\max\{C'_1,C''_1\}$ and $C_2=8\max\{C'_1,C''_2\}$,
we conclude that
\[
	\mathcal{L}(\rho_i)\leq C_1\mathcal{V}(\rho)+C_2\mathcal{D}^{d+3}(\rho)
\]
for all slices of $\rho$.
\end{proof}

We now have all the ingredients that we need to prove \thmref{mainthm}.

\begin{proof}[Proof of \thmref{mainthm}]
We first show that \condA{r} holds for all strata $H_r$. This implies,
in particular, that the second statement of \thmref{mainthm} holds for
paths starting and ending at fixed vertices. If $\rho$ is a path
starting and ending at arbitrary vertices, then \thmref{ttimproved},
\partref{fixedvert} yields that $f_\#(\rho)$ starts and ends at {\em
fixed} vertices, so that, in fact, the second statement of
\thmref{mainthm} follows from \condA{r} in this case as well.

We note that \condA{0} holds trivially, and we assume inductively that
\condA{r-1} holds for some $r$. We want to prove \condA{r}.

Assume that $H_r$ is an exponentially growing stratum, and let $\rho$ be a
smooth hallway of height $r$ such that $\rho_0$ starts and ends at fixed
vertices. If $\rho_0$ is a concatenation of Nielsen paths of height $r$
and paths in $G_{r-1}$, then we can split $\rho_0$ at the endpoints of
its subpaths in $G_{r-1}$, and \condA{r-1} completes the proof. We now
assume that $\rho_0$ is not a concatenation of Nielsen paths and paths
in $G_{r-1}$.

By \lemref{smoothlem}, we have constants $C_1,C_2$ such that
\[
	\mathcal{L}(\rho_i)\leq C_1\mathcal{V}(\rho)+C_2\mathcal{D}^{d+1}(\rho)
\]
for all slices $\rho_i$. Moreover, by \lemref{lowerbound2}, there exists
some $C>0$ and $\lambda>1$, independently of $\rho$, such that
\[
	\mathcal{V}(\rho)\geq C\lambda^{\mathcal{D}(\rho)}.
\]
We can easily find some constant $B$ such that $BC\lambda^k\geq C_2k^{d+1}$
for all $k\geq 0$. Now the first part of \condA{r} follows, with
$K_r=C_1+B$. \lemref{Cqs2} yields the second part of \condA{r}, so that
\condA{r} holds.

We now assume that $H_r$ is a polynomially growing stratum. Because of 
\propref{polyprop}, we only need to consider the following situation:
Either $H_r$ is fast, or $H_r$ is truly polynomial, but $\rho$ contains
fast polynomial edges or non-Nielsen subpaths in exponentially growing
strata.

In order to see that the second part of \condA{r} holds for a
$C$-quasismooth hallway $\rho$ of height $r$, we apply
cutting and sawtooth constructions to $\rho$, obtaining a collection of
$(C+|u_r|)$-quasismooth hallways of height $r-1$ or less, so that the
second part of \condA{r} immediately follows from the second part of
\condA{r-1}.

Now, given a smooth hallway $\rho$ of height $r$, we apply cutting and
sawtooth constructions again, obtaining a collection of
$2|u_r|$-quasismooth hallways. For each slice $\rho_i$, the second part
of \condA{r-1} yields a polynomial bound on the number of edges marked
by linear strata (\defref{markings}).  Now, since either $H_r$ is fast
or $\rho$ contains fast polynomial edges or non-Nielsen subpaths in
exponentially growing strata, \lemref{lowerbound2} provides an
exponential lower bound for the number of visible edges. As before, the
exponential lower bound for visible edges easily dominates the
polynomial lower bound for edges marked by linear strata, which
completes the proof of \condA{r}.

Finally, in order to prove the first part of \thmref{mainthm}, we need to
understand the dynamics of circuits. Let $\sigma$ be a circuit of height
$r$. If $H_r$ is a polynomially growing stratum, then \remref{howtosplit}
yields that $\sigma$ splits, at fixed vertices, into basic paths of height
$r$ and paths in $G_{r-1}$, so that \condA{r} proves the claim.

Assume that $H_r$ is an exponentially growing stratum. If $\sigma$ is
a concatenation of Nielsen paths of height $r$ and paths in $G_{r-1}$,
then we can split $\sigma$ at the endpoints of its subpaths in $G_{r-1}$,
so that \condA{r-1} completes the proof in this case. We now assume that
$\sigma$ is not a concatenation of Nielsen paths and subpaths in $G_{r-1}$.
Then $\sigma$ splits at a point $p$ in $H_r$, so that we may interpret
$\sigma$ as a path starting and ending at $v$. Let $\rho$ be a smooth
hallway with $\rho_0=\sigma$. Then, by \lemref{smoothlem}, we can find
constants $C_1,C_2$ such that
\[
	\mathcal{L}(\rho_i)\leq C_1\mathcal{V}(\rho)+C_2\mathcal{D}^d(\rho)
\]
for all slices $\rho_i$. Moreover, by \lemref{lowerbound2}, we can find
constants $C,\lambda$ such that
\[
	\mathcal{V}(\rho)\geq C\lambda^{\mathcal{D}(\rho)}.
\]
As before, we find some constant $B$ such that $BC\lambda^k\geq C_2k^{d}$
for all $k\geq 0$, so that the first statement of \thmref{mainthm} holds
with $K_r=C_1+B$.

Finally, if $\rho_0$ is a Nielsen path of height $r$, then there is
nothing to show.  This completes the proof.
\end{proof}

\bibliographystyle{alpha}
\bibliography{pb}
\par

% $Id: address.tex,v 1.3 2006/02/25 03:19:37 brinkman Exp $
{\noindent \sc
Department of Mathematics\\
The City College of CUNY\\
New York, NY 10031\\
U.S.A.\\}
\par
{\noindent \em E-mail:} brinkman@sci.ccny.cuny.edu

\end{document}